\documentclass[11pt,letterpaper]{amsart}

\usepackage{amscd,amssymb,xspace}

\usepackage{array,dcolumn}
\newcolumntype{C}{>{$}c<{$}}
\newcolumntype{L}{>{$}l<{$}}
\newcolumntype{R}{>{$}r<{$}}

\makeatletter
\let\c@equation=\c@subsection
\makeatother

\newtheorem{theorem}[subsection]{Theorem}
\newtheorem{proposition}[subsection]{Proposition}
\newtheorem{lemma}[subsection]{Lemma}
\newtheorem{corollary}[subsection]{Corollary}

\theoremstyle{definition}

\newtheorem{definition}[subsection]{Definition}

\thicklines

\newcommand{\Z}{\mathbb{Z}}

\newcommand{\Q}{\mathbb{Q}}

\newcommand{\om}{\omega}
\newcommand{\Om}{\Omega}

\renewcommand{\o}{\otimes}
\newcommand{\half}{\tfrac12}
\newcommand{\p}{\partial}

\newcommand{\bull}{\bullet}
\renewcommand{\*}{\cdot}

\renewcommand{\]}{{]\!]}}
\renewcommand{\[}{{[\![}}
\renewcommand{\)}{{)\!)}}

\newcommand{\<}{\langle}
\renewcommand{\>}{\rangle}

\DeclareMathOperator{\Edge}{\mathit{E}}

\DeclareMathOperator{\VERT}{\mathit{V}}

\newcommand{\CP}{\mathbb{CP}}
\newcommand{\EE}{\mathbb{E}}
\newcommand{\CM}{\mathcal{M}}
\newcommand{\CC}{\mathcal{C}}
\newcommand{\Mbar}{\overline{\mathcal{M}}}

\DeclareMathOperator{\Aut}{Aut}
\DeclareMathOperator{\Res}{res}
\DeclareMathOperator{\ev}{ev}
\DeclareMathOperator{\coker}{coker}
\DeclareMathOperator{\NN}{N}
\DeclareMathOperator{\NE}{NE}
\DeclareMathOperator{\ZZ}{Z}
\DeclareMathOperator{\ZE}{ZE}

\DeclareMathOperator{\supp}{supp}

\newcommand{\Dual}{\vee}
\newcommand{\Nov}{\Lambda}
\renewcommand{\SS}{S}

\newcommand{\embedding}{embedding\xspace}

\begin{document}

\title{Intersection theory on $\Mbar_{1,4}$ and elliptic Gromov-Witten
invariants}

\author{E. Getzler}

\address{Max-Planck-Institut f\"ur Mathematik, Gottfried-Claren-Str.\ 26,
D-53225 Bonn, Germany}

\curraddr{Department of Mathematics, Northwestern University, Evanston, IL
60208-2730, USA}

\email{getzler@math.nwu.edu}


\begin{abstract}
We find a new relation among codimension $2$ algebraic cycles in the moduli
space $\Mbar_{1,4}$, and use this to calculate the elliptic Gromov-Witten
invariants of projective spaces $\CP^2$ and $\CP^3$.
\end{abstract}

\maketitle

In this paper, we find a new relation among codimension $2$ algebraic
cycles in $\Mbar_{1,4}$. The main application of the new relation is to the
calculation of elliptic Gromov-Witten invariants. For example, we show that
if $V$ has no primitive cohomology in degrees above $2$, the elliptic
Gromov-Witten invariants are determined by the elliptic Gromov-Witten
invariants
$$
\<I_{1,1,\beta}^V\> : H^{2(i+1)}(V,\Q)\to\Q , \quad 0\le
i=c_1(V)\cap\beta<\dim(V) ,
$$
together with the rational Gromov-Witten invariants.

In \cite{elliptic3}, we will prove, using mixed Hodge theory, that the
cycles $[\Mbar(G)]$, as $G$ ranges over all stable graphs of genus $1$ and
valence $n$, span the even dimensional homology of $\Mbar_{1,n}$, and that
the new relation, together with those already known in genus $0$, generate
all relations among these cycles. This result is the analogue, in genus
$1$, of a theorem of Keel \cite{Keel} in genus $0$.

Our new relation is closely related to a relation in $A_2(\Mbar_3)\o\Q$
discovered by Faber (Lemma 4.4 of \cite{Faber}); the image of his relation
in $H_4(\Mbar_3,\Q)$ under the cycle map is the same as the push-forward of
our relation under the map $\Mbar_{1,4}\to\Mbar_3$ obtained by contracting
the $4$ tails pairwise. This suggests that our new relation should actually
be a rational equivalence\footnote{Since this paper was written,
Pandharipande \cite{Pandharipande} has found a direct geometric proof of
the relation of Theorem \ref{relation}, showing that it is a linear
equivalence, by means of an auxilliary moduli space of admissible covers of
$\CP^1$.}.

Let us illustrate our results with the case of the projective plane. The
genus $0$ and genus $1$ potentials of $\CP^2$ equal
\begin{align*}
F_0(\CP^2) &= \frac12 (t_0t_1^2 + t_0^2t_2) + \sum_{n=1}^\infty N^{(0)}_n q^n
e^{nt_1} \frac{t_2^{3n-1}}{(3n-1)!} , \\
F_1(\CP^2) &= - \frac{t_1}{8} + \sum_{n=1}^\infty N^{(1)}_n q^n e^{nt_1}
\frac{t_2^{3n}}{(3n)!} , \\
\end{align*}
where $t_0$, $t_1$ and $t_2$ are formal variables, of degree $-2$, $0$ and
$2$ respectively, dual to the classes $1\in H^0(\CP^2,\Q)$, $\om\in
H^2(\CP^2,\Q)$ and $\om^2\in H^4(\CP^2,\Q)$ respectively, and $N^{(0)}_n$
and $N^{(1)}_n$ are the number of rational, respectively elliptic, plane
curves of degree $n$ which meet $3n-1$, respectively $3n$, generic
points. Kontsevich and Manin \cite{KM} establish the recursion relation
$$
N^{(0)}_n = \sum_{n=i+j} \textstyle
\bigl( \binom{3n-4}{3i-2} i^2j^2 - i^3j \binom{3n-4}{3i-1} \bigr) N^{(0)}_i
N^{(0)}_j ,
$$
which, together with the initial condition $N^{(0)}_1=1$, determines the
coefficients $N^{(0)}_n$. In Section 2, we prove that the coefficients
$N^{(1)}_n$ satisfy the recursion
\begin{align} \label{recursion}
6N^{(1)}_n & = \sum_{n=i+j+k} {\textstyle
\binom{3n-2}{3j-1,3k-1} ij^3k^3 (2i-j-k) N^{(1)}_{i}
N^{(0)}_{j} N^{(0)}_{k} } \\
& {} + 2 \sum_{n=i+j} \Bigl( {\textstyle \binom{3n-2}{3i} ij^2(8i-j) -
\binom{3n-2}{3i-1} 2(i+j)j^3 } \Bigr) N^{(1)}_{i} N^{(0)}_{j} \notag \\
& {} - \frac{1}{24} \biggl( \sum_{n=i+j} {\textstyle
\binom{3n-2}{3i-1} (n^2-3n-6ij)i^3j^3 N^{(0)}_{i} N^{(0)}_{j} } + 6n^3(n-1)
N^{(0)}_n \biggr) . \notag
\end{align}
In Table 1, we list the first few coefficients $N^{(1)}_n$; for comparison,
we also include the corresponding rational Gromov-Witten invariants. We
have checked that our results for $N^{(1)}_n$ agree in degrees up to $6$
with those obtained by Caporaso and Harris \cite{CH}.

Recently, Eguchi, Hori and Ziong \cite{EHX} have proposed a bold
conjecture, generalizing the conjectured of Witten \cite{Witten} and proved
by Kontsevich \cite{KdV} that the Gromov-Witten invariants of a point (``in
the large phase space'') are the highest weight vector for a certain
Virasoro algebra. Their conjecture implies in particular the recursion
$$
N^{(1)}_n = {\textstyle\frac{1}{12} \binom{n}{3} N^{(0)}_n + \frac{1}{9}}
\sum_{n=i+j} \textstyle \binom{3n-1}{3i-1} (3i^2-2i)j N^{(0)}_i N^{(1)}_j ,
$$
which is far simpler than ours. Pandharipande \cite{Pandharipande} has
proved that this recursion is a formal consequence of \eqref{recursion}.

\begin{table} \label{CP2}
\caption{Rational and elliptic Gromov-Witten invariants of $\CP^2$}
$$\begin{tabular}{|C|R|R|} \hline
n & N^{(0)}_n \quad\quad\quad\quad & N^{(1)}_n \quad\quad\quad\quad \\  \hline
1 & 1 & 0 \\
2 & 1 & 0 \\
3 & 12 & 1 \\
4 & 620 & 225 \\
5 & 87\,304 & 87\,192 \\
6 & 26\,312\,976 & 57\,435\,240 \\
7 & 14\,616\,808\,192 & 60\,478\,511\,040 \\
8 & 13\,525\,751\,027\,392 & 96\,212\,546\,526\,096 \\ \hline
\end{tabular}$$
\end{table}

The situation for the elliptic Gromov-Witten invariants of $\CP^3$ is a
little more complicated. The genus $0$ and $1$ potentials of $\CP^3$ have
the form
\begin{align*}
F_0(\CP^3) &= \frac{t_0^2t_3}{2} + t_0t_1t_2 + \frac{t_1^3}6 +
\sum_{n=1}^\infty \sum_{4n=a+2b} N^{(0)}_{ab} q^n e^{nt_1}
\frac{t_2^at_3^b}{a!b!} , \\
F_1(\CP^3) &= - \frac{t_1}{4} + \sum_{n=1}^\infty \sum_{4n=a+2b}
N^{(1)}_{ab} q^n e^{nt_1} \frac{t_2^at_3^b}{a!b!} ,
\end{align*}
where $t_i$ is the formal variable, of degree $2i-2$, dual to $\om^i\in
H^{2i}(\CP^3,\Q)$, and $N^{(g)}_{ab}$ is the Gromov-Witten invariant which
``counts'' the stable maps of genus $g$ and degree $n$ to $\CP^3$ which
meet $a$ generic lines and $b$ generic points. As we show in Section 6, the
elliptic Gromov-Witten invariants are no longer positive integers: for
example, $N^{(1)}_{02}=-1/12$. In \cite{cp3}, we use the methods of
this paper to prove that the linear combination
$N^{(1)}_{ab}+(2n-1)N^{(0)}_{ab}/12$ counts the number of elliptic space
curves which meet $a$ generic planes and $b$ generic points.

\bigskip \noindent \textbf{Acknowledgments.} Conversations with K. Behrend,
E. Looijenga, Yu.\ Manin and especially with C. Faber, enabled me to write
this paper at all. T. Graber and R. Pandharipande informed the author of
some erroneous statements in the preprint of the paper.

I am very grateful to Yu.\ Manin, D. Zagier and the Max-Planck-Institut
f\"ur Mathematik in Bonn, where this paper was conceived, and to
A. Kupiainen and the Finnish Mathematical Society for an invitation to
Helsinki University, where much of it was finished.

The author is partially supported by the NSF.

\section{Intersection theory on $\Mbar_{1,4}$}

In this section, we calculate the relations among certain codimension two
cycles in $\Mbar_{1,4}$; one such relation was known, and we find that
there is one new one.

First, we assign names to the codimension $1$ strata of
$\Mbar_{1,4}$. Denote by $\Delta_0$ the boundary stratum of irreducible
curves in $\Mbar_{1,4}$, associated to the stable graph
$$
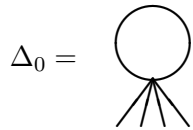
\begin{picture}(80,45)(40,745)
\put( 20,765){$\Delta_0 =$}
\put( 80,775){\circle{30}}
\put( 80,760){\line(-3,-4){ 15}}
\put( 80,760){\line(-1,-4){  5}}
\put( 80,760){\line( 1,-4){  5}}
\put( 80,760){\line( 3,-4){ 15}}
\end{picture}
$$
For each subset $S$ of $\{1,2,3,4\}$ of cardinality at least $2$, let
$\Delta_S$ be the boundary stratum associated to the stable graph with two
vertices, of genus $0$ and $1$, one edge connecting them, and with those
tails labelled by elements of $S$ attached to the vertex of genus $0$;
there are $11$ such graphs. In our pictures, we denote genus $1$ vertices
by a hollow dot, leaving genus $0$ vertices unmarked. For example,
$$
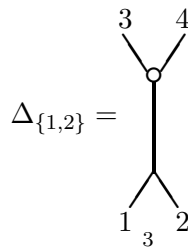
\begin{picture}(35,80)(60,722)
\put( 20,760){$\Delta_{\{1,2\}} =$}
\put( 80,780){\circle{5}}
\put( 80,777){\line( 0,-1){ 38}}
\put( 80,740){\line(-2,-3){ 10}}
\put( 80,740){\line( 2,-3){ 10}}
\put( 83,782){\line( 2, 3){ 10}}
\put( 77,782){\line(-2, 3){ 10}}
\put( 65,715){$1$}
\put( 89,715){$2$}
\put( 65,800){$3$}
\put( 89,800){$4$}
\end{picture}
$$
We only need the three $\SS_4$-invariant combinations of these $11$ strata,
which are as follows:
\begin{align*}
\Delta_2 &= \Delta_{\{1,2\}} + \Delta_{\{1,3\}} + \Delta_{\{1,4\}} +
\Delta_{\{2,3\}} + \Delta_{\{2,4\}} + \Delta_{\{3,4\}} , \\
\Delta_3 &= \Delta_{\{1,2,3\}} + \Delta_{\{1,2,4\}} + \Delta_{\{1,3,4\}} +
\Delta_{\{2,3,4\}} , \\
\Delta_4 &= \Delta_{\{1,2,3,4\}} .
\end{align*}
In summary, there are four invariant combinations of boundary strata:
$\Delta_0$, $\Delta_2$, $\Delta_3$ and $\Delta_4$.

We now turn to enumeration of the codimension two strata. These fall into
two classes, distinguished by whether they are contained in the irreducible
stratum $\Delta_0$ or not. We start by listing those which are not; each of
them is the intersection of a pair of boundary strata
$\Delta_S\*\Delta_T$. We give four examples: from these, the other strata
may be obtained by the action of $\SS_4$:
$$
\begin{picture}(100,95)(30,715)
\put(-20,755){$\Delta_{\{1,2\}} \* \Delta_{\{3,4\}} =$}
\put( 80,760){\circle{5}}
\put( 80,757){\line( 0,-1){ 18}}
\put( 80,762){\line( 0, 1){ 18}}
\put( 80,740){\line( 1,-2){ 10}}
\put( 80,740){\line(-1,-2){ 10}}
\put( 80,780){\line( 1, 2){ 10}}
\put( 80,780){\line(-1, 2){ 10}}
\put( 67,803){$1$}
\put( 87,803){$2$}
\put( 67,710){$3$}
\put( 87,710){$4$}
\end{picture}
\begin{picture}(100,85)(-20,700)
\put(-30,740){$\Delta_{\{1,2\}} \* \Delta_{\{1,2,3\}} =$}
\put( 80,760){\circle{5}}
\put( 80,757){\line( 0,-1){ 18}}
\put( 80,762){\line( 0, 1){ 18}}
\put( 80,740){\line( 1,-2){ 10}}
\put( 80,740){\line(-1,-2){ 20}}
\put( 70,720){\line( 1,-2){ 10}}
\put( 56,690){$1$}
\put( 77,690){$2$}
\put( 87,710){$3$}
\put( 77,783){$4$}
\end{picture}
$$

$$
\begin{picture}(100,75)(30,690)
\put(-40,740){$\Delta_{\{1,2\}} \* \Delta_{\{1,2,3,4\}} =$}
\put( 80,760){\circle{5}}
\put( 80,757){\line( 0,-1){ 18}}
\put( 80,740){\line( 1,-2){ 10}}
\put( 80,740){\line(-1,-2){ 20}}
\put( 70,720){\line( 1,-2){ 10}}
\put( 80,740){\line( 0,-1){ 20}}
\put( 56,690){$1$}
\put( 77,690){$2$}
\put( 77,710){$3$}
\put( 87,710){$4$}
\end{picture}
\begin{picture}(100,75)(-40,690)
\put(-50,740){$\Delta_{\{1,2,3\}} \* \Delta_{\{1,2,3,4\}} =$}
\put( 80,760){\circle{5}}
\put( 80,757){\line( 0,-1){ 18}}
\put( 80,740){\line( 1,-2){ 10}}
\put( 80,740){\line(-1,-2){ 20}}
\put( 70,720){\line( 1,-2){ 10}}
\put( 70,720){\line( 0,-1){ 20}}
\put( 56,690){$1$}
\put( 67,690){$2$}
\put( 77,690){$3$}
\put( 87,710){$4$}
\end{picture}
$$

The $\SS_4$-invariant combinations of these strata are as follows:
\begin{align*}
\Delta_{2,2} &= \Delta_{\{1,2\}} \* \Delta_{\{3,4\}} + \Delta_{\{1,3\}} \*
\Delta_{\{2,4\}} + \Delta_{\{1,4\}} \* \Delta_{\{2,3\}} , \\
\Delta_{2,3} &= \Delta_{\{1,2\}} \* \Delta_{\{1,2,3\}} + \Delta_{\{1,2\}}
\* \Delta_{\{1,2,4\}} + \Delta_{\{1,3\}} \* \Delta_{\{1,2,3\}} +
\Delta_{\{1,3\}} \* \Delta_{\{1,3,4\}} \\
& + \Delta_{\{1,4\}} \* \Delta_{\{1,2,4\}} + \Delta_{\{1,4\}} \*
\Delta_{\{1,3,4\}} + \Delta_{\{2,3\}} \* \Delta_{\{1,2,3\}} +
\Delta_{\{2,3\}} \* \Delta_{\{2,3,4\}} \\
& + \Delta_{\{2,4\}} \* \Delta_{\{1,2,4\}} + \Delta_{\{2,4\}} \*
\Delta_{\{2,3,4\}} + \Delta_{\{3,4\}} \* \Delta_{\{1,3,4\}} +
\Delta_{\{3,4\}} \* \Delta_{\{2,3,4\}} , \\
\Delta_{2,4} &= \Delta_{\{1,2\}} \* \Delta_{\{1,2,3,4\}} + \Delta_{\{1,3\}}
\* \Delta_{\{1,2,3,4\}} + \Delta_{\{1,4\}} \* \Delta_{\{1,2,3,4\}} \\
& + \Delta_{\{2,3\}} \* \Delta_{\{1,2,3,4\}} + \Delta_{\{2,4\}} \*
\Delta_{\{1,2,3,4\}} + \Delta_{\{3,4\}} \* \Delta_{\{1,2,3,4\}} , \\
\Delta_{3,4} &= \Delta_{\{1,2,3\}} \* \Delta_{\{1,2,3,4\}} +
\Delta_{\{1,2,4\}} \* \Delta_{\{1,2,3,4\}} \\
& + \Delta_{\{1,3,4\}} \* \Delta_{\{1,2,3,4\}} + \Delta_{\{2,3,4\}} \*
\Delta_{\{1,2,3,4\}} .
\end{align*}

Each of the intersections $\Delta_0\*\Delta_S$ is a codimension two stratum
in $\Delta_0$; for example
$$
\begin{picture}(80,80)(60,690)
\put(0,740){$\Delta_0\*\Delta_{\{1,2\}} =$}
\put( 80,755){\circle{30}}
\put( 80,740){\line( 1,-2){ 10}}
\put( 80,740){\line(-1,-2){ 20}}
\put( 70,720){\line( 1,-2){ 10}}
\put( 80,740){\line( 0,-1){ 20}}
\put( 56,690){$1$}
\put( 77,690){$2$}
\put( 77,710){$3$}
\put( 87,710){$4$}
\end{picture}
\begin{picture}(80,80)(0,690)
\put(-10,740){$\Delta_0\*\Delta_{\{1,2,3\}} =$}
\put( 80,755){\circle{30}}
\put( 80,740){\line( 1,-2){ 10}}
\put( 80,740){\line(-1,-2){ 20}}
\put( 70,720){\line( 1,-2){ 10}}
\put( 70,720){\line( 0,-1){ 20}}
\put( 56,690){$1$}
\put( 67,690){$2$}
\put( 77,690){$3$}
\put( 87,710){$4$}
\end{picture}
\begin{picture}(90,80)(-60,690)
\put(-15,740){$\Delta_0\*\Delta_{\{1,2,3,4\}} =$}
\put( 80,755){\circle{30}}
\put( 80,740){\line( 0,-1){ 20}}
\put( 80,720){\line(-3,-4){ 15}}
\put( 80,720){\line(-1,-4){  5}}
\put( 80,720){\line( 1,-4){  5}}
\put( 80,720){\line( 3,-4){ 15}}
\put( 61,690){$1$}
\put( 72,690){$2$}
\put( 83,690){$3$}
\put( 93,690){$4$}
\end{picture}
$$
{}From these, we may form the $\SS_4$-invariant combinations
\begin{align*}
\Delta_{0,2} &= \Delta_0\*\Delta_{\{1,2\}} + \Delta_0\*\Delta_{\{1,3\}} +
\Delta_0\*\Delta_{\{1,4\}} \\
& + \Delta_0\*\Delta_{\{2,3\}} + \Delta_0\*\Delta_{\{2,4\}} +
\Delta_0\*\Delta_{\{3,4\}} , \\
\Delta_{0,3} &= \Delta_0 \* \Delta_{\{1,2,3\}} + \Delta_0 \*
\Delta_{\{1,2,4\}} + \Delta_0 \* \Delta_{\{1,3,4\}} +
\Delta_0 \* \Delta_{\{2,3,4\}} , \\
\Delta_{0,4} &= \Delta_0 \* \Delta_{\{1,2,3,4\}} .
\end{align*}
There remain seven strata which are not expressible as intersections, which
we denote by $\Delta_{\alpha,i}$, $1\le i\le 4$, and
$\Delta_{\beta,12|34}$, $\Delta_{\beta,13|24}$ and $\Delta_{\beta,14|24}$.
We illustrate the stable graphs for two of these strata:
$$
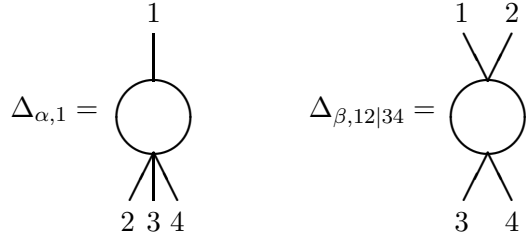
\begin{picture}(80,100)(40,720)
\put( 20,775){$\Delta_{\alpha,1} =$}
\put( 80,775){\circle{30}}
\put( 80,760){\line(-1,-2){ 10}}
\put( 80,760){\line( 0,-1){ 20}}
\put( 80,760){\line( 1,-2){ 10}}
\put( 80,790){\line( 0, 1){ 20}}
\put( 77,815){$1$}
\put( 66,727){$2$}
\put( 77,727){$3$}
\put( 87,727){$4$}
\end{picture}
\begin{picture}(100,100)(-20,720)
\put( 05,775){$\Delta_{\beta,12|34}=$}
\put( 80,775){\circle{30}}
\put( 80,760){\line(-1,-2){ 10}}
\put( 80,760){\line( 1,-2){ 10}}
\put( 80,790){\line(-1, 2){ 10}}
\put( 80,790){\line( 1, 2){ 10}}
\put( 66,815){$1$}
\put( 87,815){$2$}
\put( 66,727){$3$}
\put( 87,727){$4$}
\end{picture}
$$
Denote by $\Delta_\alpha$ and $\Delta_\beta$ the $\SS_4$-invariant
combinations of strata:
$$
\Delta_\alpha = \Delta_{\alpha,1} + \Delta_{\alpha,2} + \Delta_{\alpha,3} +
\Delta_{\alpha,4} , \quad \Delta_\beta =
\Delta_{\beta,12|34} + \Delta_{\beta,13|24} + \Delta_{\beta,14|24} .
$$

For each of these strata, let $\delta_x=[\Delta_x]$ be the corresponding
cycle in $H_\bull(\Mbar_{1,4},\Q)$, in the sense of orbifolds. (This is
sometimes denoted $[\Delta_x]_Q$ instead, but we omit the letter $Q$ from
the notation.) If the generic point of $\Delta_x$ has an automorphism group
of order $e$, then $\delta_x$ is $e^{-1}$ times the scheme-theoretic
fundamental class of $\Delta_x$; this occurs, with $e=2$, for the cycles
$\delta_{2,3}$, $\delta_{2,4}$ and $\delta_{0,4}$.
\begin{lemma} \label{trivial}
The following relation among cycles holds in $H_4(\Mbar_{1,4},\Q)$:
$$
\delta_{0,2} + 3 \delta_{0,3} + 6 \delta_{0,4} = 3 \delta_\alpha + 4
\delta_\beta .
$$
\end{lemma}
\begin{proof}
The two strata
$$
{\setlength{\unitlength}{0.01in}
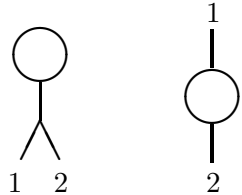
\begin{picture}(90,85)(60,685)
\put( 80,755){\circle{30}}
\put( 80,740){\line( 0,-1){ 20}}
\put( 80,720){\line(-1,-2){ 10}}
\put( 80,720){\line( 1,-2){ 10}}
\put( 63,683){$1$}
\put( 87,683){$2$}
\end{picture}
\begin{picture}(65,85)(60,727)
\put( 80,775){\circle{30}}
\put( 80,760){\line( 0,-1){ 20}}
\put( 80,790){\line( 0, 1){ 20}}
\put( 77,814){$1$}
\put( 77,725){$2$}
\end{picture}}
$$
define the same cycle, and are even rationally equivalent. (This is an
instance of the WDVV equation.) We obtain the lemma by lifting this
relation by the $6$ distinct projections $\Mbar_{1,4}\to\Mbar_{1,2}$ and
summing the answers.
\end{proof}

We can now state the main result of this section.
\begin{theorem} \label{main}
The first seven rows of the intersection matrix of the nine
$\SS_4$-in\-var\-i\-ant codimension two cycles in $\Mbar_{1,4}$ introduced
above equals
$$\begin{tabular}{C|CCCC|CCC|CC}
& \delta_{2,2} & \delta_{2,3} & \delta_{2,4} & \delta_{3,4} &
\delta_{0,2} & \delta_{0,3} & \delta_{0,4} & \delta_\alpha & \delta_\beta \\
\hline
\delta_{2,2} & 1/8 & 0 & 0 & 0 & -3 & 0 & 3/2 & 0 & 3/2 \\
\delta_{2,3} & 0 & 0 & 0 & 0 & 0 & -6 & 6 & 6 & 0 \\
\delta_{2,4} & 0 & 0 & 0 & -1/2 & 0 & 6 & -3 & 0 & 0 \\
\delta_{3,4} & 0 & 0 & -1/2 & 1/6 & 6 & -2 & 0 & 0 & 0
\\ \hline
\delta_{0,2} & -3 & 0 & 0 & 6 & 0 & 0 & 0 & 0 & 0 \\
\delta_{0,3} & 0 & -6 & 6 & -2 & 0 & 0 & 0 & 0 & 0 \\
\delta_{0,4} & 3/2 & 6 & -3 & 0 & 0 & 0 & 0 & 0 & 0 \\
\end{tabular}$$
\end{theorem}
\begin{proof}
The following lemma shows that many of the intersection numbers
vanish. (The use of this lemma simplifies our original proof of Theorem
\ref{main}, and was suggested to us by C. Faber.)
\begin{lemma}
Let $\delta$ be a cycle in $\Delta_0$. Then $\delta_0\*\delta=0$.
\end{lemma}
\begin{proof}
Consider the projection $\pi:\Mbar_{1,n}\to\Mbar_{1,1}$ which forgets all
but the first marked point, and stabilizes the marked curve which
results. The divisor $\Delta_0$ is the inverse image under $\pi$ of the
compactification divisor of $\Mbar_{1,1}$; thus, we may replace it in
calculating intersections by any cycle of the form $\pi^{-1}(x)$, where
$x\in\CM_{1,1}$. The resulting cycle has empty intersection with $\delta$,
proving the lemma.
\end{proof}

This lemma shows that all intersections among the cycles $\delta_{0,2}$,
$\delta_{0,3}$ and $\delta_{0,4}$, and between these and $\delta_\alpha$
and $\delta_\beta$ vanish.

A number of other entries in the intersection matrix vanish because the
associated strata do not meet: thus,
\begin{align*}
& \delta_{2,2}\*\delta_{2,3} = \delta_{2,2}\*\delta_{3,4} =
\delta_{2,2}\*\delta_{0,3} = \delta_{2,2}\*\delta_\alpha = 0 , \\
& \delta_{2,3}\*\delta_{2,4} = \delta_{2,3}\*\delta_\beta = 0 , \\
& \delta_{2,4}\*\delta_\alpha = \delta_{2,4}\*\delta_\beta =
\delta_{3,4}\*\delta_\alpha = \delta_{3,4}\*\delta_\beta = 0 .
\end{align*}

To calculate the remaining entries of the intersection matrix, we need the
excess intersection formula (Fulton \cite{Fulton}, Section~6.3).
\begin{proposition} \label{excess}
Let $Y$ be a smooth variety, let $X\hookrightarrow Y$ be a regular
\embedding of codimension $d$, and let $V$ be a closed subvariety of $Y$ of
dimension $n$. Suppose that the inclusion $W=X\cap V\hookrightarrow V$ is a
regular \embedding of codimension $d-e$. Then
$$
[X]\*[V] = c_e(E) \cap [W] \in A_{n-d}(W) ,
$$
where $E=(N_XY)|_W/(N_WV)$ is the \emph{excess bundle} of the intersection.
\end{proposition}

Observe that in calculating the top four rows of our intersection matrix,
at least one of the cycles which we intersect with has a regular \embedding
in $\Mbar_{1,4}$, since its dual graph is a tree. This makes the
application of the excess intersection formula straightforward.

It remains to give a formula for the normal bundles to the strata of
$\Mbar_{1,4}$.
\begin{definition}
The \emph{tautological line bundles} are defined by
$$
\om_i = \sigma_i^*\om_{\Mbar_{g,n+1}/\Mbar_{g,n}} ,
$$
where $\sigma_i:\Mbar_{g,n}\to\Mbar_{g,n+1}$, $1\le i\le n$, are the $n$
canonical sections of the universal stable curve
$\Mbar_{g,n+1}\to\Mbar_{g,n}$. Denote the Chern class $c_1(\om_i)$ by
$\psi_i$.
\end{definition}

To apply the excess intersection formula, we need to know the normal
bundles of strata $\Mbar(G)\subset\Mbar_{g,n}$. The following result gives
a partial answer to this question, and is all that we need for the
calculations in this paper: a proof may be found in Section~4 of
Hain-Looijenga \cite{HL}.
\begin{proposition} \label{normal}
Let $G$ be a stable graph of genus $g$ and valence $n$, and let
$$
p: \prod_{v\in\VERT(G)} \Mbar_{g(v),n(v)} \to \Mbar_{g,n}
$$
be the ramified cover (of degree $|\Aut(G)|$) of the closed stratum
$\Mbar(G)$ of $\CM_{g,n}$. Each edge $e$ of the graph determines two flags
$s(e)$ and $t(e)$, and hence two tautological line bundles $\om_{s(e)}$ and
$\om_{t(e)}$ on $\prod_{v\in\VERT(G)} \Mbar_{g(v),n(v)}$, and the normal
bundle of $p$ is given by the formula
$$
N_p = \bigoplus_{e\in\Edge(G)} \om_{s(e)}^\Dual\o\om_{t(e)}^\Dual .
\qed$$
\end{proposition}

In particular, if the graph $G$ has no automorphisms, so that $p$ is an
\embedding, the bundle $N_p$ may be identified with the normal bundle of
the stratum $\Mbar(G)$.

It is now straightforward to calculate the remaining entries of the
intersection matrix. We will use the integrals
\begin{equation} \label{tau}
\int_{\Mbar_{0,4}} \psi_i = 1 , \quad \int_{\Mbar_{1,1}} \psi_1 =
\int_{\Mbar_{1,2}} \psi_1 \cup \psi_2 = \frac{1}{24} ,
\end{equation}
which are proved in Witten \cite{Witten}.

In performing the calculations, it is helpful to introduce a graphical
notation for the cycle obtained from a stratum by capping with a monomial
in the Chern classes $-\psi_i$: we point a small arrow along each flag $i$
where we intersect by the class $-\psi_i$. (This notation generalizes that
of Kaufmann \cite{Kaufmann}, who considers the case of trees where the
genus of each vertex is $0$. The minus signs come from the inversion
accompanying the tautological line bundles in the formula of Proposition
\ref{normal}.) One then calculates the contribution of such a graph by
multiplying together factors for each vertex equal to the integral over
$\Mbar_{g(v),n(v)}$ of the appropriate monomial in the classes $-\psi_i$,
and dividing by the order of the automorphism group $\Aut(G)$: in
particular, this vanishes unless there are $3(g(v)-1)+n(v)$ arrows at each
vertex $v$.

We illustrate the sort of enumeration which arises with one of the most
complicated of these calculations, that of
$\delta_{2,4}\*\delta_{2,4}$. Two sorts of terms contribute: $6$~terms of
the form
$$
\bigl( \delta_{\{1,2\}}\*\delta_{\{1,2,3,4\}} \bigr)^2 = \frac{1}{24} ,
$$
and $6$ terms of the form
$$
\delta_{\{1,2\}}\*\delta_{\{1,2,3,4\}}\*\delta_{\{3,4\}}\*\delta_{\{1,2,3,4\}}
= - \frac{1}{24} .
$$
Applying the excess intersection formula, we see that
$$
\bigl(\delta_{\{1,2\}}\*\delta_{\{1,2,3,4\}}\bigr)^2 =
c_2\bigl(N_{\Delta_{\{1,2\}}\cap\Delta_{\{1,2,3,4\}}}\Mbar_{1,4}\bigr) \cap
\bigl( \delta_{\{1,2\}}\*\delta_{\{1,2,3,4\}} \bigr) .
$$
Expanding the second Chern class of the normal bundle, we see that each
term contributes the sum of four graphs:
$$
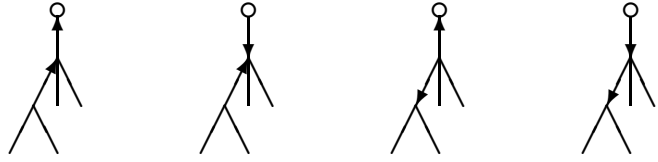
\begin{picture}(80,90)(60,690)
\put( 80,760){\circle{5}}
\put( 80,740){\vector( 0, 1){ 18}}
\put( 80,740){\line( 1,-2){ 10}}
\put( 80,740){\line( 0,-1){ 20}}
\put( 70,720){\vector( 1, 2){ 10}}
\put( 70,720){\line(-1,-2){ 10}}
\put( 70,720){\line( 1,-2){ 10}}
\end{picture}
\begin{picture}(80,90)(60,690)
\put( 80,760){\circle{5}}
\put( 80,757){\vector( 0,-1){ 18}}
\put( 80,740){\line( 1,-2){ 10}}
\put( 80,740){\line( 0,-1){ 20}}
\put( 70,720){\vector( 1, 2){ 10}}
\put( 70,720){\line(-1,-2){ 10}}
\put( 70,720){\line( 1,-2){ 10}}
\end{picture}
\begin{picture}(80,90)(60,690)
\put( 80,760){\circle{5}}
\put( 80,740){\vector( 0, 1){ 18}}
\put( 80,740){\line( 1,-2){ 10}}
\put( 80,740){\line( 0,-1){ 20}}
\put( 80,740){\vector(-1,-2){ 10}}
\put( 70,720){\line(-1,-2){ 10}}
\put( 70,720){\line( 1,-2){ 10}}
\end{picture}
\begin{picture}(50,90)(60,690)
\put( 80,760){\circle{5}}
\put( 80,757){\vector( 0,-1){ 18}}
\put( 80,740){\line( 1,-2){ 10}}
\put( 80,740){\line( 0,-1){ 20}}
\put( 80,740){\vector(-1,-2){ 10}}
\put( 70,720){\line(-1,-2){ 10}}
\put( 70,720){\line( 1,-2){ 10}}
\end{picture}
$$
Only the first graph is nonzero, since in the other cases, the wrong number
of arrows point towards the vertices. And the first graph contributes
$$
\int_{\Mbar_{0,4}} (-\psi_1) \* \int_{\Mbar_{1,1}} (-\psi_1) = \frac{1}{24} .
$$

In the case of terms of the form
$\delta_{\{1,2\}}\*\delta_{\{1,2,3,4\}}
\*\delta_{\{3,4\}}\*\delta_{\{1,2,3,4\}}$,
the excess dimension $e$ equals $1$, and we must calculate the degree of
the excess bundle on the stratum
$\Delta_{\{1,2\}}\cap\Delta_{\{3,4\}}\cap\Delta_{\{1,2,3,4\}}$. Two graphs
contribute:
$$
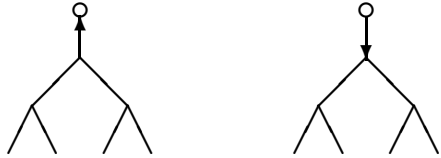
\begin{picture}(120,90)(60,690)
\put( 80,760){\circle{5}}
\put( 80,740){\vector( 0, 1){ 18}}
\put( 80,740){\line( 1,-1){ 20}}
\put( 80,740){\line(-1,-1){ 20}}
\put( 60,720){\line(-1,-2){ 10}}
\put( 60,720){\line( 1,-2){ 10}}
\put(100,720){\line(-1,-2){ 10}}
\put(100,720){\line( 1,-2){ 10}}
\end{picture}
\begin{picture}(60,90)(60,690)
\put( 80,760){\circle{5}}
\put( 80,757){\vector( 0,-1){ 18}}
\put( 80,740){\line( 1,-1){ 20}}
\put( 80,740){\line(-1,-1){ 20}}
\put( 60,720){\line(-1,-2){ 10}}
\put( 60,720){\line( 1,-2){ 10}}
\put(100,720){\line(-1,-2){ 10}}
\put(100,720){\line( 1,-2){ 10}}
\end{picture}
$$
Only the first of these graphs gives a nonzero value, namely
$$
\int_{\Mbar_{1,1}} (-\psi_1) = - \frac{1}{24} .
$$
This completes our outline of the proof of Theorem \ref{main}.
\end{proof}

The intersection matrix of Theorem \ref{main} has rank $7$. We now apply
the results of \cite{genus1}, where we calculated the character of the
$\SS_n$-modules $H^i(\Mbar_{1,n},\Q)$: these calculations show that $\dim
H^4(\Mbar_{1,4},\Q)^{\SS_4}=7$. This shows that our $9$ cycles span
$H^4(\Mbar_{1,4},\Q)^{\SS_4}$, and that the nullspace of the intersection
matrix gives relations among them. We already know one such relation, by
Lemma \ref{trivial}. Calculating the remaining null-vector of the
intersection matrix, we obtain the main theorem of this paper.
\begin{theorem} \label{relation}
The following new relation among cycles holds:
$$
12\delta_{2,2} - 4\delta_{2,3} - 2\delta_{2,4} + 6\delta_{3,4} +
\delta_{0,3} + \delta_{0,4} - 2\delta_\beta = 0 .
\qed$$
\end{theorem}

Using this theorem, it is easy to calculate the remaining intersections
among our $9$ strata:
$$
\delta_\alpha\*\delta_\alpha=16 , \quad \delta_\alpha\*\delta_\beta=-12 ,
\quad \delta_\beta\*\delta_\beta=9 .
$$
C. Faber informs us that the direct calculation of these intersection
numbers is not difficult. This would allow a different approach to the
proof of Theorem \ref{relation}, using the theorem of \cite{elliptic3} that
the strata of $\Mbar_{1,n}$ span the even-dimensional rational cohomology.

\section{Gromov-Witten invariants}

In the remainder of this paper, we apply the new relation to the
calculation of elliptic Gromov-Witten invariants: we will do this
explicitly for curves and for the projective plane $\CP^2$, and prove some
general results in other cases.

\subsection{The Novikov ring}
Let $V$ be a smooth projective variety of dimension $d$. In studying the
Gromov-Witten invariants, it is convenient to work with cohomology with
coefficients in the Novikov ring $\Nov$ of $V$, which we now define.

Let $\NN_1(V)$ be the abelian group
$$
\NN_1(V) = \ZZ_1(V) / \text{numerical equivalence} ,
$$
and let $\NE_1(V)$ be its sub-semigroup
$$
\NE_1(V) = \ZE_1(V) / \text{numerical equivalence} ,
$$
where $\ZZ_1(V)$ is the abelian group of $1$-cycles on $V$, and $\ZE_1(V)$
is the semigroup of effective $1$-cycles. (Recall that two $1$-cycles $x$
and $y$ are numerically equivalent $x\equiv y$ when $x\*Z=y\*Z$ for any
Cartier divisor $Z$ on $V$.)

The Novikov ring is
\begin{align*}
\Nov &= \Q[\NN_1(V)] \o_{\Q[\NE_1(V)]} \Q\[\NE_1(V)\] \\
&= \textstyle \bigl\{ a = \sum_{\beta\in\NN_1(V)} a_\beta q^\beta \mid
\text{ $\supp(a) \subset \beta_0+\NE_1(V)$ for some  $\beta_0\in\NN_1(V)$}
\bigr\} ,
\end{align*}
with product $q^{\beta_1}q^{\beta_2}=q^{\beta_1+\beta_2}$ and grading
$|q^\beta|=-2c_1(V)\cap\beta$. That the product is well-defined is shown by
the following proposition (Koll\'ar \cite{Kollar}, Proposition II.4.8).
\begin{proposition} \label{Mori}
If $V$ is a projective variety with K\"ahler form $\om$, the set
$$
\{\beta\in\NE_1(V)\mid \om\cap\beta\le c\}
$$
is finite for each $c>0$.
\qed\end{proposition}

For example, if $V=\CP^n$, then $\NN_1(\CP^n)=\Z\*[L]$, where $[L]$ is the
cycle defined by a line $L\subset\CP^n$, and $\Nov\cong\Q\(q\)$, with
grading $|q|=-2(n+1)$, since $c_1(\CP^n)\cap[L]=n+1$.

If $V=E$ is an elliptic curve, then $\NN_1(E)=\Z\*[E]$, and
$\Nov\cong\Q\(q\)$, concentrated in degree $0$.

\subsection{Stable maps}
The definition of Gromov-Witten invariants is based on the study of the
moduli stacks $\Mbar_{g,n}(V,\beta)$ of stable maps of Kontsevich, which
have been shown by Behrend and Manin \cite{BM} to be complete
Deligne-Mumford stacks (though not in general smooth).

For each $N\ge0$, let $\pi_{n,N}:\Mbar_{g,n+N}(V,\beta) \to
\Mbar_{g,n}(V,\beta)$ be the projection which forgets the last $N$ marked
points of the stable curve, and stabilizes the resulting map.

In the
special case $N=1$, we obtain a fibration
$$
\pi: \Mbar_{g,n+1}(V,\beta) \to \Mbar_{g,n}(V,\beta)
$$
which is shown by Behrend and Manin to be the universal curve; that is, its
fibre over a stable map $(f:C\to V,x_i)$ is the curve $C$. Denote by
$f:\Mbar_{g,n+1}(V,\beta)\to V$ the universal stable map, obtained by
evaluation at $x_{n+1}$.

\subsection{The virtual fundamental class}
There are projections $\Mbar_{g,n}(V,\beta)\to\Mbar_{g,n}$, when
$2(g-1)+n>0$, which send the stable map $(f:C\to V,x_i)$ to the
stabilization of $(C,x_i)$. If the sheaf $R^1\pi_*f^*TV$ vanishes on
$\Mbar_{g,n}(V,\beta)$, the Riemann-Roch theorem predicts that the fibres
of the projection $\Mbar_{g,n}(V,\beta)\to\Mbar_{g,n}$ have dimension
$$
d(1-g)+c_1(V)\cap\beta ,
$$
and hence that $\Mbar_{g,n}(V,\beta)$ has dimension
$$
d(1-g)+c_1(V)\cap\beta + \dim\Mbar_{g,n} = (3-d)(1-g)+c_1(V)\cap\beta + n .
$$
This hypothesis is only rarely true, and in any case only in genus
$0$. However, Behrend-Fantecchi \cite{B,BF} and Li-Tian \cite{LT} show that
there is a bivariant class
$$
[\Mbar_{g,n}(V,\beta)/\Mbar_{g,n},R^\bull\pi_*f^*TV] \in
A^{d(1-g)+c_1(V)\cap\beta}(\Mbar_{g,n}(V,\beta)\to\Mbar_{g,n}) ,
$$
the virtual relative fundamental class, which stands in for
$[\Mbar_{g,n}(V,\beta)/\Mbar_{g,n}]$ in the obstructed case.

The following result is proved in \cite{B}, and sometimes permits the
explicit calculation of Gromov-Witten invariants, as we will see later.
\begin{proposition} \label{excess-virtual}
If the coherent sheaf $R^1\pi_*f^*TV$ on $\Mbar_{g,n}(V,\beta)$ is locally
trivial of dimension $e$ (the \emph{excess dimension}), then
$\Mbar_{g,n}(V,\beta)$ is smooth of dimension
$$
(3-d)(1-g)+c_1(V)\cap\beta+n+e ,
$$
and $[\Mbar_{g,n}(V,\beta)/\Mbar_{g,n},R^\bull\pi_*f^*TV] =
c_e(R^1\pi_*f^*TV) \cap [\Mbar_{g,n}(V,\beta)/\Mbar_{g,n}]$.
\qed
\end{proposition}

\subsection{Gromov-Witten invariants}
The Gromov-Witten invariant of genus $g\ge0$, valence $n\ge0$ and degree
$\beta\in\NE_1(V)$ is a cohomology operation
$$
I_{g,n,\beta}^V : H^{2d(1-g)+2c_1(V)\cap\beta+\bull}(V^n,\Q) \to
H^\bull(\Mbar_{g,n},\Q) ,
$$
defined by the formula
$$
I_{g,n,\beta}^V(\alpha_1,\dots,\alpha_n) =
[\Mbar_{g,n}(V,\beta)/\Mbar_{g,n},R^\bull\pi_*f^*TV] \cap
\ev^*(\alpha_1\boxtimes\dots\boxtimes\alpha_n) ,
$$
where $\ev:\Mbar_{g,n}(V,\beta)\to V^n$ is evaluation at the marked points:
$$
\ev : (f:C\to V,x_i) \mapsto (f(x_1),\dots,f(x_n)) \in V^n .
$$

Capping $I_{g,n,\beta}^V$ with the fundamental class $[\Mbar_{g,n}]$, we
obtain a numerical invariant
$$
\<I_{g,n,\beta}^V\> : H^{2(d-3)(1-g)+2c_1(V)\cap\beta+2n}(V^n,\Q) \to \Q .
$$
This is the $n$-point correlation function of two-dimensional topological
gravity with the topological $\sigma$-model associated to $V$ as a
background \cite{Witten}. Note that if $\beta\ne0$, $\<I_{g,n,\beta}^V\>$
may be defined even when $2(g-1)+n\le0$, even though $I_{g,n,\beta}^V$ does
not exist.

Introducing the Novikov ring, we may define the generating function
$$
I_{g,n}^V = \sum_{\beta\in\NE_1(V)} q^\beta I_{g,n,\beta}^V :
H^*(V,\Nov)^{\o n} \to H^\bull(\Mbar_{g,n},\Nov) ,
$$
along with its integral over the fundamental class $[\Mbar_{g,n}]$
$$
\<I_{g,n}^V\> = \sum_{\beta\in\NE_1(V)} q^\beta \<I_{g,n,\beta}^V\> :
H^*(V,\Nov)^{\o n} \to \Nov ,
$$
Note that $I_{g,n}^V$ and $\<I_{g,n}^V\>$ are invariant under the action of
the symmetric group $\SS_n$ on $H^*(V,\Nov)^{\o n}$.

In the special case of zero degree, the moduli space $\Mbar_{g,n}(V,\beta)$
is isomorphic to $\Mbar_{g,n}\times V$. This allows us to calculate the
Gromov-Witten invariants $\<I_{0,3,0}^V\>$ and $\<I_{1,1,0}^V\>$. The
former is given by the explicit formula
$$
\<I_{0,3,0}^V(\alpha_1,\alpha_2,\alpha_3)\> = \int_V
\alpha_1\cup\alpha_2\cup\alpha_2 .
$$
This formula is very simple to prove, since the moduli space
$\Mbar_{0,3}(V,0)\cong V$ is smooth, with dimension equal to its virtual
dimension $d$, and thus the virtual fundamental class
$[\Mbar_{0,3}(V,0),R^\bull\pi_*f^*TV]$ may be identified with the
fundamental class of $V$. A similar proof shows that $\<I_{0,n,0}^V\>$
vanishes if $n>3$.

The calculation of the Gromov-Witten invariant $\<I_{1,1,0}^V\>$ (see
Bershadsky et al.\ \cite{BCOV}) is a good illustration of the application
of Proposition \ref{excess-virtual}.
\begin{proposition} \label{BCOV}
$$
\<I_{1,1,0}^V(\alpha)\> = -\frac{1}{24} \int_V c_{d-1}(V)\cup\alpha ,
$$
while $\<I_{1,n,0}^V\>=0$ if $n>1$.
\end{proposition}
\begin{proof}
The moduli stack $\Mbar_{1,n}(V,0)$ is isomorphic to $\Mbar_{1,n}\times V$,
and the obstruction bundle $R^1\pi_*f^*TV$ is isomorphic to the vector
bundle $\EE^\Dual\boxtimes TV$, of rank $d$, where
$\EE=\pi_*\om_{\Mbar_{1,n+1}/\Mbar_{1,n}}$. Hence $R^1\pi_*f^*TV$ has top
Chern class
$$
c_d(\EE^\Dual\o f^*TV) = 1\boxtimes f^*c_d(V) - \lambda_1 \boxtimes
f^*c_{d-1}(V) ,
$$
where $\lambda_1=c_1(\EE)$. By Proposition \ref{excess-virtual},
\begin{align*}
\<I_{1,n,0}^V(\alpha_1,\dots,\alpha_n)\> &= \int_{\Mbar_{1,n}\times V}
c_d(\EE^\Dual\o f^*TV) \boxtimes (\alpha_1\cup\dots\cup\alpha_n) \\
&= - \int_{\Mbar_{1,n}} \lambda_1 \* \int_V c_{d-1}(V) \cup
\alpha_1\cup\dots\cup\alpha_n .
\end{align*}
On dimensional grounds, $\<I_{1,n,0}^V\>$ vanishes if $n>1$, while the
formula follows when $n=1$ from $\lambda_1\cap[\Mbar_{1,1}]=\frac{1}{24}$.
\end{proof}

\subsection{The puncture axiom}
One of the basic axioms satisfied by Gromov-Witten invariants is expressed
in the relationship between virtual fundamental classes
$$
[\Mbar_{g,n+1}(V,\beta)/\Mbar_{g,n+1},R^\bull\pi_*f^*TV] = \pi^*
[\Mbar_{g,n}(V,\beta)/\Mbar_{g,n},R^\bull\pi_*f^*TV] .
$$
Here, $\pi^*:A^k(\Mbar_{g,n}(V,\beta)\to\Mbar_{g,n}) \to
A^k(\Mbar_{g,n+1}(V,\beta)\to\Mbar_{g,n+1})$ is the operation of flat
pullback associated to the diagram
$$\begin{CD}
\Mbar_{g,n+1}(V,\beta) @>>> \Mbar_{g,n+1} \\
@V{\pi}VV @V{\pi}VV \\
\Mbar_{g,n}(V,\beta) @>>> \Mbar_{g,n}
\end{CD}$$
This axiom implies that if $\alpha$ is a cohomology class on $V$ of degree
at most $2$ and $2(g-1)+n>0$,
\begin{equation} \label{low}
I_{g,n+1,\beta}^V(\alpha,\alpha_1,\dots,\alpha_n) =
\begin{cases}
0 , & |\alpha|=0,1 , \\
(\alpha\cap\beta) I_{g,n,\beta}^V(\alpha_1,\dots,\alpha_n) , &
|\alpha|=2 .
\end{cases}
\end{equation}

\subsection{Generating functions}
Let $\Lambda\[H\]$ be the power series ring
$\Lambda\[H_{\bull+2}(V,\Q)\]$. Let $\{\gamma^a\}_{a=0}^k$ be a homogeneous
basis of the graded vector space $H^\bull(V,\Q)$, with $\gamma^0=1$, and
let $\{t_a\}_{a=0}^k$ be the dual basis; the (homological) degree of $t_a$
equals the (cohomological) degree of $\gamma^a$ minus $2$. We may identify
the ring $\Lambda\[H\]$ with $\Lambda\[t_0,\dots,t_k\]$.

Let $F_g(V)$ be the generating function
$$
F_g(V) = \sum_{n=0}^\infty \<I_{g,n}^V\> \in \Lambda\[H\] .
$$
This is a power series of degree $2(d-3)(1-g)$. This suggests assigning to
Planck's constant $\hbar$ the degree $2(d-3)(g-1)$, and forming the total
generating function, homogeneous of degree $0$,
$$
F(V) = \sum_{g=0}^\infty \hbar^{g-1} F_g(V) .
$$

\subsection{The composition axiom}
The composition axiom for Gromov-Witten invariants gives a formula for the
integral of the Gromov-Witten invariant $I_{g,n}^V$ over the cycle
$[\Mbar(G)]$ associated to a stable graph $G$ which bears a strong
resemblance to the Feynman rules of quantum field theory:

Let $\eta_{ab}$ be the Poincar\'e form of $V$ with respect to the basis
$\{\gamma^a\}_{a=0}^k$ of $H^\bull(V,\Q)$. Then
$$
\int_{\Mbar(G)} I_{g,n}^V(\alpha_1,\dots,\alpha_n) = \frac{1}{\Aut(G)}
\sum_{\substack{a(e),b(e)=0 \\ e\in\Edge(G)}}^k \prod_{e\in\Edge(G)}
\eta_{a(e),b(e)} \prod_{v\in\VERT(G)} \<I_{g(v),n(v)}^V(\dots)\> .
$$
Here, the Gromov-Witten invariant $\<I_{g(v),n(v)}^V(\dots)\>$ is evaluated
on the cohomology classes $\alpha_i$ corresponding to the tails of $G$
which meet the vertex $v$, on the $\gamma^{a(e)}$ corresponding to edges
$e$ which start at the vertex $v$, and on the $\gamma^{b(e)}$ corresponding
to edges $e$ which end at $v$. (The right-hand side is independent of the
chosen orientation of the edges, by the symmetry of the Poincar\'e form.)

\subsection{Relations among Gromov-Witten invariants}
Let $G$ be a stable graph of genus $g$ and valence $n$. The subvariety
$\pi_{n,N}^{-1}\bigl(\Mbar(G)\bigr) \subset \Mbar_{g,n+N}$ is the union of
strata associated to the set of stable graphs obtained from $G$ by
adjoining $N$ tails $\{n+1,\dots,n+N\}$ in all possible ways to the
vertices of $G$.

For example, consider the stratum $\Delta_{12|34}\subset\Mbar_{0,4}$,
associated to the stable graph
$$
\begin{picture}(35,85)(60,715)
\put( 20,760){$\Delta_{12|34} =$}
\put( 80,775){\line( 0,-1){ 30}}
\put( 80,745){\line(-2,-3){ 10}}
\put( 80,745){\line( 2,-3){ 10}}
\put( 80,775){\line( 2, 3){ 10}}
\put( 80,775){\line(-2, 3){ 10}}
\put( 65,715){$1$}
\put( 89,715){$2$}
\put( 65,795){$3$}
\put( 89,795){$4$}
\end{picture}
$$
The inverse image $\pi_{4,N}^{-1}(\Delta_{12|34})$ consists of the union of
all strata in $\Mbar_{0,4+N}$ associated to stable graphs
$$
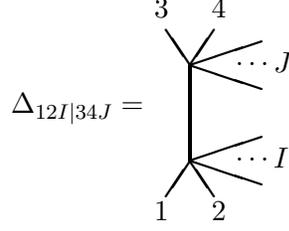
\begin{picture}(35,85)(60,715)
\put(  5,760){$\Delta_{12I|34J} =$}
\put( 80,780){\line( 0,-1){ 40}}
\put( 80,740){\line(-2,-3){ 10}}
\put( 80,740){\line( 2,-3){ 10}}
\put( 80,780){\line( 2, 3){ 10}}
\put( 80,780){\line(-2, 3){ 10}}
\put( 65,715){$1$}
\put( 89,715){$2$}
\put( 65,800){$3$}
\put( 89,800){$4$}
\put( 80,780){\line( 3,-1){ 30}}
\put( 80,780){\line( 3, 1){ 30}}
\put( 95,780){\dots}
\put(115,777){$J$}
\put( 80,740){\line( 3,-1){ 30}}
\put( 80,740){\line( 3, 1){ 30}}
\put( 95,740){\dots}
\put(115,737){$I$}
\end{picture}
$$
where $I$ and $J$ form a partition of the set $\{5,\dots,N+4\}$.

If $\delta$ is a cycle in $\Mbar_{g,n}$, define the generating function
$$
F(\delta,V) = \sum_{N=0}^\infty \int_{\pi^{-1}(\delta)} I_{g,n+N}^V :
H^{\bull+2}(V,\Nov)^{\o n} \to \Nov\[H\] .
$$
More explicitly,
\begin{multline*}
F(\delta,V)(\alpha_1,\dots,\alpha_n) \\ = \sum_{N=0}^\infty \frac{1}{N!}
\sum_{a_1,\dots,a_N} t_{a_N}\dots t_{a_1} \int_\delta \bigl( \pi_{n,N}
\bigr)_*
I_{g,n+N}^V(\gamma^{a_1},\dots,\gamma^{a_N},\alpha_1,\dots,\alpha_n) .
\end{multline*}
In particular, if $\delta=[\Mbar(G)]$ where $G$ is a stable graph, we set
$$
F(G,V)=F([\Mbar(G)],V) .
$$
If $g>1$, $F_g(V)$ is a special case of this construction, with
$\delta=[\Mbar_{g,0}]$.

A little exercise involving Leibniz's rule shows that the composition axiom
implies the following formula for these generatings functions:
\begin{equation} \label{composition}
F(G,V) = \frac{1}{\Aut(G)} \sum_{\substack{a(e),b(e)=0 \\ e\in\Edge(G)}}^k
\prod_{e\in\Edge(G)} \eta_{a(e),b(e)} \prod_{v\in\VERT(G)}
\p^{n(v)}F_{g(v)}(V) (\dots) ,
\end{equation}
where as before, the multilinear form $\p^{n(v)}F_{g(v)}(V)$ is evaluated
on the cohomology classes $\alpha_i$ corresponding to the tails of $G$
meeting the vertex $v$, on the $\gamma^{a(e)}$ corresponding to edges $e$
which start at the vertex $v$, and on the $\gamma^{b(e)}$ corresponding to
edges $e$ which end at $v$.

The composition axiom implies that any relation among the cycles
$[\Mbar(G)]$ is reflected in a relation among Gromov-Witten invariants,
which, by \eqref{composition} may be translated into a differential
equation among generating functions $F_g(V)$. An example is the rational
equivalence of the cycles associated to the three strata of $\Mbar_{0,4}$
of codimension $1$:
$$
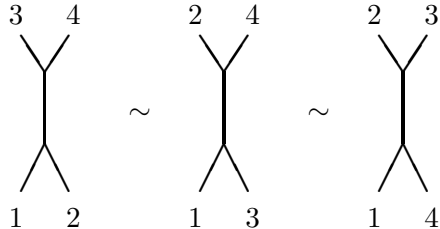
\begin{picture}(35,85)(60,715)
\put( 80,775){\line( 0,-1){ 30}}
\put( 80,745){\line(-1,-2){ 10}}
\put( 80,745){\line( 1,-2){ 10}}
\put( 80,775){\line( 2, 3){ 10}}
\put( 80,775){\line(-2, 3){ 10}}
\put( 65,710){$1$}
\put( 89,710){$2$}
\put( 65,795){$3$}
\put( 89,795){$4$}
\end{picture}
\hskip0.5in
\begin{picture}(35,85)(60,715)
\put( 40,755){$\sim$}
\put( 80,775){\line( 0,-1){ 30}}
\put( 80,745){\line(-1,-2){ 10}}
\put( 80,745){\line( 1,-2){ 10}}
\put( 80,775){\line( 2, 3){ 10}}
\put( 80,775){\line(-2, 3){ 10}}
\put( 65,710){$1$}
\put( 89,710){$3$}
\put( 65,795){$2$}
\put( 89,795){$4$}
\end{picture}
\hskip0.5in
\begin{picture}(35,85)(60,715)
\put( 40,755){$\sim$}
\put( 80,775){\line( 0,-1){ 30}}
\put( 80,745){\line(-1,-2){ 10}}
\put( 80,745){\line( 1,-2){ 10}}
\put( 80,775){\line( 2, 3){ 10}}
\put( 80,775){\line(-2, 3){ 10}}
\put( 65,710){$1$}
\put( 89,710){$4$}
\put( 65,795){$2$}
\put( 89,795){$3$}
\end{picture}
$$
The equality of the Gromov-Witten invariant $F(\delta,V)$ when evaluated on
these cycles is the WDVV equation.

In order to express the relation among the Gromov-Witten invariants implied
by Theorem \ref{relation}, it is useful to introduce certain operators
which act on elements of $\Nov[H]\o\Nov\[H\]$ through differentiation in
the first factor: the Laplacian
$$
\Delta = \frac12 \sum_{a,b=0}^k \eta_{ab} \frac{\p^2}{\p t_a\p t_b} ,
$$
and the sequence of bilinear differential operators $\Gamma_n$ by
$\Gamma_0(f,g)=fg$ and
$$
\Gamma_n(f,g) = \frac{1}{n} \bigl( \Delta\Gamma_{n-1}(f,g) -
\Gamma_{n-1}(\Delta f,g) - \Gamma_{n-1}(f,\Delta g) \bigr) .
$$
(We will abbreviate $\Gamma_1(f,g)$ to $\Gamma(f,g)$.)

\begin{proposition} \label{Relation}
Denote the derivative $\p^{n(v)}F_{g(v)}(V)/n(v)!\in\Nov[H]\o\Nov\[H\]$ by
$f_{g,n}$. (Note that $f_{g,n}=F([\Mbar_{g,n}],V)$.) Then
\begin{align*}
6 \, \Gamma(\Gamma_1(f_{1,2},f_{0,3}),f_{0,3})
&- 5 \, \Gamma(f_{1,2},\Gamma(f_{0,3},f_{0,3})) \\
& {} - 2 \, \Gamma(f_{0,3},\Gamma(f_{1,1},f_{0,4}))
+ 6 \, \Gamma(f_{0,4},\Gamma(f_{1,1},f_{0,3})) \\
& {} + \Gamma(f_{0,4},\Delta f_{0,4})
+ \Gamma(f_{0,5},\Delta f_{0,3})
- \Gamma_2(f_{0,4},f_{0,4}) = 0 .
\end{align*}
\end{proposition}
\begin{proof}
This follows from the following table, which is obtained by application of
\eqref{composition}.
$$\begin{tabular}{|L|L||L|L|}
\hline
\delta & F(\delta,V) & & \\ \hline
\delta_{2,2} & \frac12 \Gamma(\Gamma(f_{1,2},f_{0,3}),f_{0,3}) &
\delta_{0,2} & \Gamma(f_{0,3},\Delta f_{0,5}) \\
& {} - \frac14 \Gamma(f_{1,2},\Gamma(f_{0,3},f_{0,3})) &
\delta_{0,3} & \Gamma(f_{0,4},\Delta f_{0,4}) \\
\delta_{2,3} & \frac12 \Gamma(f_{1,2},\Gamma(f_{0,3},f_{0,3})) &
\delta_{0,4} & \Gamma(f_{0,5},\Delta f_{0,3}) \\
\delta_{2,4} & \Gamma(f_{0,3},\Gamma(f_{1,1},f_{0,4})) &
\delta_\alpha & \Gamma_2(f_{0,3},f_{0,5}) \\
\delta_{3,4} & \Gamma(f_{0,4},\Gamma(f_{1,1},f_{0,3})) &
\delta_\beta & \frac12 \Gamma_2(f_{0,4},f_{0,4}) \\[1pt] \hline
\end{tabular}$$
\end{proof}

When we apply Proposition \ref{Relation} with $V=\CP^2$ and evaluate the
resulting multilinear form to $\om^{\boxtimes4}$, we obtain the recursion
relation \eqref{recursion} for the elliptic Gromov-Witten invariants
$N^{(1)}_n$ of $\CP^2$.

\section{The symbol of the new relation}

We may introduce a filtration on Gromov-Witten invariants with respect to
which the leading order of our new relation takes a relatively simple form;
by analogy with the case of differential operators, we call this leading
order relation the symbol of the full relation. In some cases, this symbol
may be used to prove that elliptic Gromov-Witten invariants are determined
by rational ones.

\begin{definition}
The \emph{symbol} of a relation $\delta=0$ among cycles of strata in
$\Mbar_{g,n}$ is the set of relations among Gromov-Witten invariants
obtained by taking, for each $\beta\in\NE_1(V)$, the coefficient of
$q^\beta$ in $I_{g,n}^V\cap[\delta]$, expanding in Feynman diagrams using
the composition axiom, and setting all Gromov-Witten invariants
$\<I_{g',n',\beta'}^V\>$ other than $\<I_{g,n,\beta}^V\>$ and
$\<I_{0,3,0}^V\>$ to zero.
\end{definition}

We define a total order on the symbols $\<I_{g,n,\beta}^V\>$ by setting
$\<I_{g',n',\beta'}^V\>\prec\<I_{g,n,\beta}^V\>$ if $g'<g$, or $g'=g$ and
$n'<n$, or $g'=g$, $n'=n$ and $\beta=\beta'+\beta''$ where
$\beta''\in\NE_1(V)$ is non-zero. Thus, knowledge of the symbol determines
relations among Gromov-Witten invariants such that the error in the
relation on $\<I_{g,n,\beta}^V\>$ involves invariants
$\<I_{g',n',\beta'}^V\>$ with
$\<I_{g',n',\beta'}^V\>\prec\<I_{g,n,\beta}^V\>$. (Here, we must of course
exclude $\<I_{0,3,0}^V\>$.) We use the symbol $\sim$ to denote this
equivalence relation.

For example, the symbol of the WDVV equation is
$$
(a,b,c\cup d) + (a\cup b,c,d) \sim (-1)^{|a|(|b|+|c|)} \bigl( (b,c,a\cup d)
+ (b\cup c,a,d) \bigr) ,
$$
where we have abbreviated $\<I_{0,n,\beta}^V
(\alpha_1,\alpha_2,\alpha_3,\alpha_4,\alpha_5,\dots,\alpha_n)\>$ to
$(\alpha_1,\alpha_2,\alpha_3,\alpha_4)$.

Next, consider the symbol of the relation
$$
\pi_{4,n-4}^{-1}\bigl(12\delta_{2,2} - 4\delta_{2,3} - 2\delta_{2,4} +
6\delta_{3,4} + \delta_{0,3} + \delta_{0,4} - 2\beta\bigr) = 0
$$
in $H_{2n-4}(\Mbar_{1,n},\Q)$. Only the cycles $\delta_{2,2}$ and
$\delta_{2,3}$ contribute terms to the symbol. Abbreviate the Gromov-Witten
class $\<I_{1,n,\beta}^V(\alpha_1,\alpha_2,\alpha_3,\dots,\alpha_n)\>$ to
$\{\alpha_1,\alpha_2\}$. Up to a numerical factor to be determined, the
cycle $\delta_{2,2}$ contributes the expression
$$
\{a\cup b,c\cup d\} + (-1)^{|b|\,|c|} \{a\cup c,b\cup d\} +
(-1)^{(|b|+|c|)|d|} \{a\cup d,b\cup c\} .
$$
This numerical factor equals
$$
\frac{1}{24} \* 3 \* 12 \* 8 = 12 .
$$
The factor $1/24$ comes from symmetrization over the four inputs, the
factor of $3$ from the three strata making up $\delta_{2,2}$, the factor of
$12$ is the coefficient of the cycle in the relation, and the factor $8$ is
illustrated by listing all of the graphs which contribute a term $\{a\cup b
,c\cup d\}$:
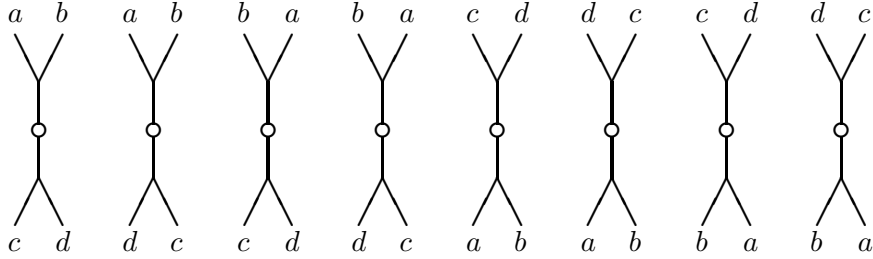
$$
\def\nnn{\begin{picture}(48,95)(56,715)
\put( 80,760){\circle{5}}
\put( 80,757){\line( 0,-1){ 18}}
\put( 80,762){\line( 0, 1){ 18}}
\put( 80,740){\line( 1,-2){ 10}}
\put( 80,740){\line(-1,-2){ 10}}
\put( 80,780){\line( 1, 2){ 10}}
\put( 80,780){\line(-1, 2){ 10}}}
\nnn
\put( 67,805){$a$}
\put( 87,805){$b$}
\put( 67,709){$c$}
\put( 87,709){$d$}
\end{picture}
\nnn
\put( 67,805){$a$}
\put( 87,805){$b$}
\put( 67,709){$d$}
\put( 87,709){$c$}
\end{picture}
\nnn
\put( 67,805){$b$}
\put( 87,805){$a$}
\put( 67,709){$c$}
\put( 87,709){$d$}
\end{picture}
\nnn
\put( 67,805){$b$}
\put( 87,805){$a$}
\put( 67,709){$d$}
\put( 87,709){$c$}
\end{picture}
\nnn
\put( 67,805){$c$}
\put( 87,805){$d$}
\put( 67,709){$a$}
\put( 87,709){$b$}
\end{picture}
\nnn
\put( 67,805){$d$}
\put( 87,805){$c$}
\put( 67,709){$a$}
\put( 87,709){$b$}
\end{picture}
\nnn
\put( 67,805){$c$}
\put( 87,805){$d$}
\put( 67,709){$b$}
\put( 87,709){$a$}
\end{picture}
\nnn
\put( 67,805){$d$}
\put( 87,805){$c$}
\put( 67,709){$b$}
\put( 87,709){$a$}
\end{picture}
$$
Similarly, the cycle $\delta_{2,3}$ contributes the expression
\begin{multline*}
\{a,b\cup c\cup d\} + (-1)^{|a|\,|b|} \{b,a\cup c\cup d\} \\
+ (-1)^{(|a|+|b|)|c|} \{c,a\cup b\cup d\} + (-1)^{(|a|+|b|+|c|)|d|}
\{d,a\cup b\cup c\} ,
\end{multline*}
with numerical factor
$$
\frac{1}{24} \* 12 \* (-4) \* 6 = - 12 ;
$$
the factor $12$ counts the strata making up $\delta_{2,3}$, $-4$ is the
coefficient of the cycle in the relation, and we illustrate the factor $6$
by listing all of the graphs which contribute a term $\{a,b\cup c\cup d\}$:
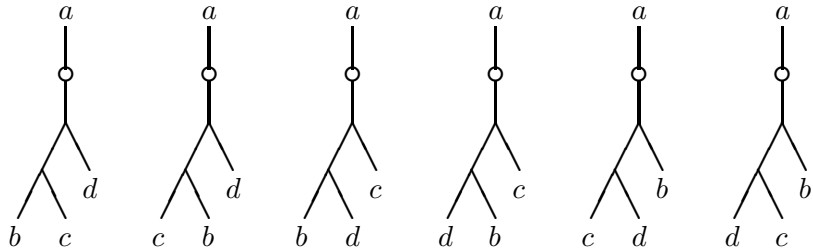
$$
\def\nnn{\begin{picture}(60,95)(50,695)
\put( 80,760){\circle{5}}
\put( 80,757){\line( 0,-1){ 18}}
\put( 80,762){\line( 0, 1){ 18}}
\put( 80,740){\line( 1,-2){ 10}}
\put( 80,740){\line(-1,-2){ 20}}
\put( 70,720){\line( 1,-2){ 10}}
}
\nnn
\put( 56,688){$b$}
\put( 77,688){$c$}
\put( 87,708){$d$}
\put( 77,783){$a$}
\end{picture}
\nnn
\put( 56,688){$c$}
\put( 77,688){$b$}
\put( 87,708){$d$}
\put( 77,783){$a$}
\end{picture}
\nnn
\put( 56,688){$b$}
\put( 77,688){$d$}
\put( 87,708){$c$}
\put( 77,783){$a$}
\end{picture}
\nnn
\put( 56,688){$d$}
\put( 77,688){$b$}
\put( 87,708){$c$}
\put( 77,783){$a$}
\end{picture}
\nnn
\put( 56,688){$c$}
\put( 77,688){$d$}
\put( 87,708){$b$}
\put( 77,783){$a$}
\end{picture}
\nnn
\put( 56,688){$d$}
\put( 77,688){$c$}
\put( 87,708){$b$}
\put( 77,783){$a$}
\end{picture}
$$
In conclusion, we obtain the following result.
\begin{theorem} \label{symbol}
Abbreviating
$\<I_{1,n,\beta}^V(\alpha_1,\alpha_2,\alpha_3,\dots,\alpha_n)\>$ to
$\{\alpha_1,\alpha_2\}$, we have
\begin{align*}
\Psi(a,b,c,d) & = \{a\cup b,c\cup d\} + (-1)^{|b|\,|c|} \{a\cup c,b\cup d\} +
(-1)^{(|b|+|c|)|d|} \{a\cup d,b\cup c\} \\
& {} - \{a,b\cup c\cup d\} - (-1)^{|a|\,|b|} \{b,a\cup c\cup d\} \\
& {} - (-1)^{(|a|+|b|)|c|} \{c,a\cup b\cup d\} - (-1)^{(|a|+|b|+|c|)|d|}
\{d,a\cup b\cup c\} \sim 0 .
\end{align*}
\end{theorem}

Note that the linear form $\Psi(a,b,c,d)$ is (graded) symmetric in its four
arguments, and vanishes if any of them equals $1$.

\begin{corollary} \label{reduce}
If $\om\in H^2(V,\Q)$ and $a,b\in H^\bull(V,\Q)$, then for $j>2$,
$$\textstyle
\{\om^i\cup a,\om^{j-i}\cup b\} = \binom{i+2}{2} \{a,\om^j\cup b\} .
$$
\end{corollary}
\begin{proof}
By Theorem \ref{symbol}, we have for $i\ge0$ and $j>2$,
\begin{multline*}
\Psi(\om,\om^{i+1}\cup a,\om,\om^{j-i-3}\cup b) -
\Psi(\om,\om^i\cup a,\om,\om^{j-i-2}\cup b) \\
\begin{aligned} \sim & \quad \bigl(
2\{\om^{i+2}\cup a,\om^{j-i-2}\cup b\} +
\{\om^2,\om^{j-2}\cup a\cup b\} \\
{} & {} - \{\om^{i+1}\cup a,\om^{j-i-1}\cup b\} -
\{\om^{i+3}\cup a,\om^{j-i-3}\cup b\} \bigr) \\
{} & {} - \bigl( 2\{\om^{i+1}\cup a,\om^{j-i-1}\cup b\} +
\{\om^2,\om^{j-2}\cup a\cup b\} \\
{} & {} - \{\om^i\cup a,\om^{j-i}\cup b\}
- \{\om^{i+2}\cup a,\om^{j-i-2}\cup b\} \bigr) \\
\sim & \quad \{\om^i\cup a,\om^{j-i}\cup b\} -
3\{\om^{i+1}\cup a,\om^{j-i-1}\cup b\} \\
{} & \quad {} + 3\{\om^{i+2}\cup a,\om^{j-i-2}\cup b\}
- \{\om^{i+3}\cup a,\om^{j-i-3}\cup b\} \sim 0 .
\end{aligned}
\end{multline*}
This implies that the function
$a(i,j)=\{\om^i\cup a,\om^{j-i}\cup b\}$ satisfies the difference
equation
$$
a(i,j) - 3 a(i+1,j) + 3a(i+2,j) - a(i+3,j) \sim 0
$$
with solution $a(i,j)\sim\binom{i+2}{2} a(0,j)$.
\end{proof}

We can now prove a weak analogue for elliptic Gromov-Witten invariants of
the (first) Reconstruction Theorem of Kontsevich-Manin (Theorem 3.1 of
\cite{KM}). For $0\le j\le d$, let $P_j(V)=\coker( H^{j-2}(V,\Q)
\xrightarrow{\om\cup\*} H^j(V,\Q) )$ be the $j^{\text{th}}$ primitive
cohomology group of $V$.
\begin{theorem} \label{reconstruction}
If $P^i(V)=0$ for $i>2$, the elliptic Gromov-Witten invariants of $V$ are
determined by its rational Gromov-Witten invariants together with the
Gromov-Witten invariants $\<I_{1,1,\beta}(-)\>:H^{2i+2}(V,\Q)\to\Q$ for
$0\le c_1(V)\cap\beta=i<d$. (These are all of the non-vanishing
Gromov-Witten invariants $\<I_{1,1,\beta}(\alpha)\>$.)
\end{theorem}
\begin{proof}
We proceed by induction: by hypothesis, $\<I_{g,n,\beta}^V\>$ is known for
$g=0$ or $g=1$ and $n=1$. Now consider the Gromov-Witten invariant
$\<I_{1,n,\beta}^V(\alpha_1,\dots,\alpha_n)\>$, where $n>1$. By
\eqref{low}, we may assume that $|\alpha_i|>2$, and under the hypotheses of
the proposition, we may write it as $\om^{p_i}\cup\gamma_i$ where
$|\gamma_i|\le2$ is a primitive cohomology class.

Step 1: If any two indices $p_i$ and $p_j$ satisfy $p_i+p_j>2$, we may
apply Corollary \ref{reduce} to replace the pair
$(\om^{p_i}\cup\gamma_i,\om^{p_j}\cup\gamma_j)$ by
$(\gamma_i,\om^{p_i+p_j}\cup\gamma_j)$. If $|\gamma_1|=1$, the result
vanishes by \eqref{low}, while if $|\gamma_1|=2$, we may apply \eqref{low}
to reduce $n$ by $1$.

Step 2: We are reduced to considering
$\<I_{1,n,\beta}^V(\om\cup\gamma_1,\dots,\om\cup\gamma_n)\>$, where the
classes $\gamma_i$ have degree $1$ or degree $2$. Applying Theorem
\ref{symbol}, we see that
$$
\Psi(\om,\gamma_1,\om,\gamma_2) = 2\{\om\cup\gamma_1,\om\cup\gamma_2\} +
\{\om^2,\gamma_1\cup\gamma_2\} \sim 0 .
$$
In particular, we may assume that $n=2$, since otherwise, we would be able
to return to Step 1. There are two cases.

Step 2a: If the classes $\gamma_i$ are both of degree $1$, we see that
$\{\om\cup\gamma_1,\om\cup\gamma_2\}\sim0$, since in that case
$\gamma_1\cup\gamma_2$ has degree $2$ and we may apply \eqref{low}.

Step 2b: If the classes $\gamma_i$ are both of degree $2$, there is a class
$\gamma\in H^2(V,\Q)$ such that $\gamma_1\cup\gamma_2 = \om\cup\gamma$,
since $P^4(V)=0$. We must calculate
$$
\<I_{1,2,\beta}^V(\om^2,\gamma_1\cup\gamma_2)\>
= \<I_{1,2,\beta}^V(\om^2,\om\cup\gamma)\> \sim 6
\<I_{1,2,\beta}^V(1,\om^3\cup\gamma)\> = 0 ,
$$
where we have applied Corollary \ref{reduce} and \eqref{low}.
\end{proof}

Two special cases of this result are worth singling out:
\begin{enumerate}
\item If $V$ is a surface, the elliptic Gromov-Witten invariants are
determined by the rational invariants together with
$\<I_{1,1,\beta}(-)\>:H^2(V,\Q)\to\Q$ for $c_1(V)\cap\beta=0$ and
$\<I_{1,1,\beta}(-)\>:H^4(V,\Q)\to\Q$ for $c_1(V)\cap\beta=1$. If $V$ is
the blow-up of $\CP^2$ at a finite number of points, only $\beta=0$
satisfies $c_1(V)\cap\beta<2$, and by Proposition \ref{BCOV},
$\<I_{1,1,0}\>$ is determined by $c_1(V)$, while the rational Gromov-Witten
invariants are determined by the WDVV equation (G\"ottsche-Pandharipande
\cite{GP}).
\item If $V=\CP^d$, the elliptic Gromov-Witten invariants are determined by
the rational Gromov-Witten invariants.
\end{enumerate}

\section{Gromov-Witten invariants of curves}

To illustrate our new relation, we start with the case where $V$ is a
curve. We will only discuss curves of genus $0$ and $1$, since for curves
of higher genus, $I_{g,n,\beta}^V=0$ if $\beta\ne0$, and the new relation
is identically satisfied.

\subsection{The projective line}
When $V=\CP^1$, the potential $F_g$ is a power series of degree $4g-4$ in
variables $t_0$ and $t_1$ (of degree $-2$ and $0$) and the generator $q$ of
$\Lambda$, of degree $-4=-2c_1(\CP^1)\cap[\CP^1]$. By degree counting,
together with \eqref{low}, we see that
$$
F_g(\CP^1) =
\begin{cases} \displaystyle
t_0^2t_1/2 + q e^{t_1} , & g=0 , \\
\displaystyle -t_1/24 , & g=1 , \\ 0 , & g>1 ;
\end{cases}$$
the only thing which is not immediate is the coefficient of $q$ in
$F_0(\CP^1)$, which is the number of maps of degree $1$ from $\CP^1$ to
itself, up to isomorphism, and clearly equals $1$.

It is easy to calculate $F(\delta,\CP^1)$ for $\delta$ equal to one of our
nine $2$-cycles: all of them vanish except
$$
F(\delta_{3,4},\CP^1) = \frac{t_1^4}{24} \o(-qe^{t_1}/6) ;
F(\delta_{0,4},\CP^1) = \frac{t_1^4}{24} \o qe^{t_1} ;
F(\delta_\alpha,\CP^1) = \frac{t_1^4}{24} \o 2qe^{t_1} .
$$
We see that the new relation holds among these potentials.

\subsection{Elliptic curves}
Let $E$ be an elliptic curve. Denote by $\xi,\eta$ variables of degree $-1$
corresponding to a basis of $H_1(E,\Z)$ such that $\<\xi,\eta\>=1$. The
ring $\Lambda$ has one generator $q$, of degree $0$ (since
$c_1(V)=0$). Since there are no rational curves in $E$ of positive degree,
we have
$$
F_0(E) = t_0^2t_1/2 + t_0\xi\eta .
$$
It is shown in \cite{BCOV} that
\begin{equation} \label{Eisenstein}
F_1(E) = - \frac{t_1}{24} + \sum_{\beta=1}^\infty \frac{\sigma(\beta)}{\beta}
q^\beta \bigl(e^{\beta t_1} - 1\bigr) ,
\end{equation}
since $\<I_{1,1,\beta}^E(\om)\>=\sigma(\beta)$ counts the number of
unramified covers of degree $\beta$ of the curve $E$ up to automorphisms,
which are easily enumerated. An equivalent form of \eqref{Eisenstein} is
$$
\frac{\p F_1(E)}{\p t_1} = G_2(qe^{t_1}) ,
$$
where
$$
G_2(q) = - \frac{1}{24} + \sum_{\beta=1}^\infty \sigma(\beta) q^\beta
$$
is the Eisenstein series of weight $2$. By degree counting, we also see
that $F_g(E)=0$ for $g>1$.

Note that the Gromov-Witten invariants of an elliptic curve are invariant
under deformation; this is true for any smooth projective variety $V$
(Li-Tian \cite{LT}). In fact, the definition of Gromov-Witten invariants
extends to any almost-K\"ahler manifold (a symplectic manifold with
compatible almost-complex structure), and the resulting invariants are
independent of the almost-complex structure (Li-Tian \cite{LT:symp}).

It is simple to calculate the Gromov-Witten potentials $F(\delta,E)$ for
our nine $2$-cycles in $\Mbar_{1,4}$.
\begin{lemma} \label{elliptic}
We have
$$
F(\delta_{2,2},E) = \Bigl( \frac{5}{12} G_4(qe^{t_1}) - G_2(qe^{t_1})^2
\Bigr) (t_0t_1+\xi\eta)^2 = \frac{q}{2} (t_0t_1+\xi\eta)^2 + O(q^2) ,
$$
$F(\delta_{2,3},E)=3F(\delta_{2,2},E)$, while the remaining $7$ potentials
vanish.
\qed\end{lemma}

Again, we see that the new relation holds.

\section{The Gromov-Witten invariants of $\CP^2$}

The Gromov-Witten potential $F_g(\CP^2)$ is a power series of degree $2g-2$
in variables $t_0$, $t_1$ and $t_2$, of degrees $-2$, $0$ and $2$, where
$t_i$ is dual to $\om^i$, and the generator $q$ of $\Lambda$, of degree
$-6=-2c_1(\CP^2)\cap[L]$.

By degree counting, together with \eqref{low}, we see that
$$
F_g(\CP^2) =
\begin{cases}
\displaystyle \frac{1}{2} (t_0t_1^2 + t_0^2t_2) + \sum_{\beta=1}^\infty
N^{(0)}_\beta q^\beta e^{\beta t_1} \frac{t_2^{3\beta-1}}{(3\beta-1)!} , &
g=0 , \\
\displaystyle - \frac{t_1}{8} + \sum_{\beta=1}^\infty N^{(1)}_\beta q^\beta
e^{\beta t_1} \frac{t_2^{3\beta}}{(3\beta)!} , & g=1 , \\
\displaystyle \sum_{\beta=1}^\infty N^{(g)}_\beta q^\beta e^{\beta t_1}
\frac{t_2^{3\beta+g-1}}{(3\beta+g-1)!} , & g>1 , \end{cases}
$$
where $N^{(g)}_\beta$ are certain rational coefficients.

Using the Severi theory of plane curves, we will show that $N^{(g)}_\beta$
is the answer to an enumerative problem for plane curves; in particular, it
is a non-negative integer. This phenomenon is special to del~Pezzo
surfaces: we have already seen that the elliptic Gromov-Witten invariants
of an elliptic curve are non-integral, while for $\CP^3$, they are not even
positive.

We apply the following result, which is Proposition 2.2 of Harris
\cite{Harris}.
\begin{proposition} \label{Harris}
Let $S$ be a smooth rational surface. Let $\pi:\CC\to\CM$ be a family of
curves of geometric genus $g$ with $\CM$ irreducible, and let $f:\CC\to\CM$
be a map such that on each component of a general fibre $\CC_z$ of $\pi$,
the restriction $f_z$ of $f$ to $\CC_z$ is not constant and $f_z^*\om_S$
has negative degree.

Let $W$ be the image of the map from $\CM$ to the Chow variety of curves on
$S$ defined by sending $z\in\CM$ to the curve $\CC_z$. Then
$\dim(W)\le-\deg(f_z^*\om_S)+g-1$, and if equality holds, then $f_z$ is
birational for all $z\in\CM$.
\qed
\end{proposition}

\begin{corollary} \label{Severi}
The coefficient $N^{(g)}_\beta$ equals the number of irreducible plane
curves of arithmetic genus $g$ and degree $\beta$ passing through
$3\beta+g-1$ general points in $\CP^2$.
\end{corollary}
\begin{proof}
Let $\CM$ be a component of the boundary
$\Mbar_{g,n}(\CP^2,\beta)\setminus\CM_{g,n}(\CP^2,\beta)$, and consider the
family of curves $\CC\to\CM$ obtained by restricting the universal curve
$\Mbar_{g+1,n}(\CP^2,\beta)\to\Mbar_{g,n}(\CP^2,\beta)$ to $\CM$ and
contracting to a point all components of the fibres on which $f$ has degree
$0$.

The geometric genus of the fibres of this family is bounded above by $g-1$.
Applying Proposition \ref{Harris}, we see that the image of $\CM$ in the Chow
variety of plane curves has dimension at most $3\beta+g-2$.

On the other hand, if $\CM$ is a component of $\CM_{g,n}(\CP^2,\beta)$, and
$\CC\to\CM$ is the universal family of curves $\CC\to\CM$, we see that the
image of $\CM$ in the Chow variety of plane curves has dimension less than
$3\beta+g-1$ unless the stable maps parametrized by $\CM$ are birational to
their image.

The Gromov-Witten invariant $N^{(g)}_\beta$ equals the degree of the
intersection of the image of $\Mbar_{g,3\beta+g-1}(\CP^2,\beta)$ in the
Chow variety of curves in $\CP^2$ with the cycle of curves passing through
$3\beta+g-1$ general points. By Bertini's theorem for homogenous spaces
\cite{Kleiman}, we see that the points of intersection are reduced and lie
in the components of $\CM_{g,n}(\CP^2,\beta)$ on which the map $f$ is
birational to its image, and hence an \embedding. (This argument is
borrowed from Section 6 of Fulton-Pandharipande \cite{FP}.) The result
follows.
\end{proof}

\subsection{Comparison with the formulas of Caporaso and Harris}
Caporaso and Harris \cite{CH} have calculated the numbers $N^{(g)}_\beta$
for all $g\ge0$, and we now turn the comparison of our results for
$N^{(1)}_\beta$ . We have not been able to find a proof that our answers
agree, but we have verified that this is so for $\beta\le6$.

The recursion of Caporaso and Harris for the Gromov-Witten invariants of
$\CP^2$ is more easily applied if it is recast in terms of generating
functions.

\begin{definition}
If $\alpha$ is a partition, denote by $\ell(\alpha)$ the number of parts of
$\alpha$ and by $|\alpha|$ the sum $\alpha_1+\dots+\alpha_{\ell(\alpha)}$
of the parts of $\alpha$. Let $\alpha!$ be the product $\alpha! = \alpha_1!
\dots \alpha_{\ell(\alpha)}!$.
\end{definition}

Fix a line $L$ in $\CP^2$. If $\alpha$ and $\beta$ are partitions with
$|\alpha|+|\beta|=d$, and $\Om$ is a collection of $\ell(\alpha)$ general
points of $L$, let
$V^{d,\delta}(\alpha,\beta)(\Om)=V^{d,\delta}(\alpha,\beta)$ be the
generalized Severi variety: the closure of the locus of reduced plane
curves of degree $d$ not containing $L$, smooth except for $\delta$ double
points, having order of contact $\alpha_i$ with $L$ at $\Om_i$, and to
order $\beta_1,\dots,\beta_{\ell(\beta)}$ at $\ell(\beta)$ further
unassigned points of $L$. For example, $V^{d,\delta}(0,1^d)$ is the
classical Severi variety of plane curves of degree $d$ with $\delta$ double
points, while $V^{d,\delta}(0,21^{d-1})$ is the closure of the locus of
plane curves tangent to $L$ at a smooth point.

Denote by $V_0^{d,\delta}(\alpha,\beta)$ the union of the components of
$V^{d,\delta}(\alpha,\beta)$ whose general point is an irreducible
curve. Let $N^{d,\delta}(\alpha,\beta)$ be the degree of
$V^{d,\delta}(\alpha,\beta)$ and let $N_0^{d,\delta}(\alpha,\beta)$ be the
degree of $V_0^{d,\delta}(\alpha,\beta)$. Form the generating functions
\begin{align*}
Z &= \sum \frac{z^{\binom{d+1}{2}-\delta+\ell(\beta)}}
{\bigl(\binom{d+1}{2}-\delta+\ell(\beta)\bigr)!} \frac{p^\alpha}{\alpha!}
q^\beta N^{d,\delta}(\alpha,\beta) , \\
F &= \sum \frac{z^{\binom{d+1}{2}-\delta+\ell(\beta)}}
{\bigl(\binom{d+1}{2}-\delta+\ell(\beta)\bigr)!} \frac{p^\alpha}{\alpha!}
q^\beta N_0^{d,\delta}(\alpha,\beta) .
\end{align*}
The integer $\binom{d+1}{2}-\delta+\ell(\beta)$ is the dimension of the
variety $V^{d,\delta}(\alpha,\beta)$. The union of curves of degree $d_i$,
$1\le i\le n$, with $\delta_i$ double points and partitions $\alpha_i$ and
$\beta_i$ is a (reducible) curve has degree $d=d_1+\dots+d_n$ with
$$
\delta = \delta_1 + \dots + \delta_n + \sum_{i<j} \delta_i\delta_j
$$
double points and partitions $\alpha=(\alpha_1,\dots,\alpha_n)$ and
$\beta=(\beta_1,\dots,\beta_n)$. This formula for $\delta$ amounts to the
condition that the sum of the dimensions of the generalized Severi
varieties $V_0^{d_i,\delta_i}(\alpha_i,\beta_i)$ equals the dimension of
$V^{d,\delta}(\alpha,\beta)$. The proof of the relationship $Z = \exp(F)$
between these two generating functions is an exercise in the definition of
degree (see Ran \cite{Ran}).

Caporaso and Harris prove a recursion which in terms of the generating
function $Z$ may be written
$$
\frac{\p Z}{\p z} = \sum_{k=1}^\infty kq_k\frac{\p Z}{\p p_k} + \Res_{t=0}
\biggl[ \exp\Bigl( \sum_{k=1}^\infty t^{-k} p_k + \sum_{k=1}^\infty k t^k
\frac{\p}{\p q_k} \Bigr) \biggr] Z ,
$$
where $\Res_{t=0}$ is the residue with respect to the formal variable $t$,
in other words, the coefficient of $t^{-1}$ when the exponential is
expanded.%
\footnote{The resemblance of the right-hand side to the Hamiltonian of the
Liouville model is striking --- we have no idea why operators so closely
resembling vertex operators make their appearance here.} Dividing by $Z$,
we obtain
$$
\frac{\p F}{\p z} = \sum_{k=1}^\infty kq_k\frac{\p F}{\p p_k} + \Res_{t=0}
\biggl[ \exp\Bigl( \sum_{k=1}^\infty t^{-k} p_k + F|_{q_k\mapsto q_k+kt^k}
- F \Bigr) \biggr] ,
$$
which clearly allows the recursive calculation of the coefficients
$N_0^{d,\delta}(\alpha,\beta)$.

As a special case of $Z=\exp(F)$, we have
$$
1 + \sum \frac{z^{\binom{d+2}{2}-\delta-1} q^d N^{d,\delta}}
{\bigl(\binom{d+2}{2}-\delta-1\bigr)!} = \exp \biggl( \sum
\frac{z^{\binom{d+2}{2}-\delta-1} q^d N_0^{d,\delta}}
{\bigl(\binom{d+2}{2}-\delta-1\bigr)!} \biggr) ,
$$
since $\binom{d+1}{2}-\delta+d=\binom{d+2}{2}-\delta-1$. Expanding the
exponential, we obtain
$$
N^{d,\delta} = \sum_{n=1}^\infty \frac{1}{n!} \sum_{d=d_1+\dots+d_n} \\
\sum_{\substack{\delta=\sum_{i<j}\delta_i\delta_j\\+\delta_1+\dots+\delta_n}}
\frac{\bigl(\binom{d+2}{2}-\delta-1\bigr)! N_0^{d_1,\delta_1} \dots
N_0^{d_n,\delta_n}} {\bigl(\binom{d_1+2}{2}-\delta_1-1\bigr)! \dots
\bigl(\binom{d_i+2}{2}-\delta_i-1\bigr)!} .
$$
For example, with $d=5$, we obtain
\begin{align*}
N_0^{5,4} &= N^{5,4} - \frac{16!}{14!2!} N_0^{4,0}N_0^{1,0} = 36975 - 120 \*
1 = 36855 , \\
N_0^{5,5} &= N^{5,5} - \frac{15!}{13!2!} N_0^{4,1}N_0^{1,0} = 90027 - 105
\* 27 = 87192 ,
\end{align*}
while with $d=6$ and $\delta=9$, we obtain
\begin{align*}
N_0^{6,9} &= N^{6,9} - 18! \biggl( \frac{N_0^{5,4}}{16!}
\frac{N_0^{1,0}}{2!} - \frac{N_0^{4,1}}{13!} \frac{N_0^{2,0}}{3!}
- \frac{1}{2} \Bigl( \frac{N_0^{3,0}}{9!} \Bigr)^2 - \frac{1}{2}
\frac{N_0^{4,0}}{14!} \Bigl( \frac{N_0^{1,0}}{2!} \Bigr)^2 \biggr) \\
&= 63338881 - 153 \* 36855 \* 1 + 8568 \* 27 \* 1 + \half \* 48620 \* 1^2 +
\half \* 18360 \* 1 \* 1^2 \\ &= 57435240
\end{align*}
in agreement with the recursion \eqref{recursion}.

By Proposition \ref{Severi}, the relation between the numbers
$N_0^{d,\delta}$ and the Gromov-Witten invariants is very simple:
$N^{(g)}_d=N_0^{d,\delta}$ where $g=\binom{d-1}{2}-\delta$. In terms of
$F$, the Gromov-Witten potentials $F_g(\CP^2)$ are given by the formula
$$
\sum_{g=0}^\infty \hbar^{g-1} F_g(\CP^2) = \frac{1}{2\hbar}
(t_0^2t_2+t_0t_1^2) - \frac{t_1}{8} +
F\big|_{\substack{(q_1,q_2,\dots)=(\hbar^{-3}qe^{t_1},0,\dots) \\
(p_1,p_2,\dots)=(0,0,\dots) , z=\hbar t_2}} .
$$

\section{The elliptic Gromov-Witten invariants of $\CP^3$}

For $g=0$ and $g=1$, the Gromov-Witten potentials of the projective space
$\CP^3$ have the form
$$
F_g(\CP^3) = \begin{cases}
\bigl( \frac{1}{2} t_0^2t_3 + t_0t_1t_2 + \frac{1}{6} t_1^3 \bigr) +
\displaystyle \sum_{4\beta=a+2b} N^{(0)}_{ab} q^\beta e^{\beta t_1}
\frac{t_2^at_3^b}{a!b!} , & g=0 , \\
\displaystyle - \frac{t_1}{4} + \sum_{4\beta=a+2b}
N^{(1)}_{ab} q^\beta e^{\beta t_1} \frac{t_2^at_3^b}{a!b!} , & g=1 .
\end{cases}$$
Here, $t_i$ is the formal variable of degree $2i-2$ dual to $\om^i\in
H^{2i}(\CP^3,\Q)$ and $q$ is the generator of the Novikov ring
$\Lambda\cong\Q\(q\)$ of $\CP^3$. By Proposition \ref{BCOV}, the
coefficient of $t_1$ in $F_1(\CP^3)$ equals $-c_2(\CP^3)/24$.

Thus, the coefficient $N^{(g)}_{ab}$ is a rational number which ``counts''
the number of stable maps of degree $\beta$ from a curve of genus $g$ to
$\CP^3$ meeting $a$ generic lines and $b$ generic points.

It is shown by Fulton and Pandharipande \cite{FP} that $N^{(0)}_{ab}$
equals the number of rational space curves of degree $\beta$ which meet $a$
generic lines and $b$ generic points. In particular, they are non-negative
integers. By contrast, the coefficients $N^{(1)}_{ab}$ are neither positive
nor integral: for example, $N^{(1)}_{02}=-1/12$. In \cite{cp3}, we prove
the following result.
\begin{theorem}
The number of elliptic space curves of degree $\beta$ passing through $a$
generic lines and $b$ generic points, where $4\beta=a+2b$, equals
$N^{(1)}_{ab} + (2\beta-1)N^{(0)}_{ab}/12$.
\end{theorem}

By evaluating the equation of Proposition \ref{Relation} on
$\om\boxtimes\om\boxtimes\om\boxtimes\om$, we obtain the following relation
among the elliptic Gromov-Witten for $\CP^3$: if $a\ge2$, then
\begin{multline*}
3 N^{(1)}_{ab} = 4 nN^{(1)}_{a-2,b+1} - \tfrac{1}{4} n^2 N^{(0)}_{ab} +
\tfrac{1}{6} n^3 (n-3) N^{(0)}_{a-2,b+1} \\
\shoveleft{ {} - 2 \sum_{\substack{a-2=a_1+a_2\\b+1=b_1+b_2}}
\textstyle N^{(1)}_{a_1b_1} N^{(0)}_{a_2b_2} n_2^2 (n-3n_1)
\binom{a-2}{a_1} \Bigr\{ n_1 \binom{b}{b_1} + n_2 \binom{b}{b_1-1} \Bigr\} } \\
\shoveleft{ {} + \sum_{\substack{a=a_1+a_2\\b=b_1+b_2}} N^{(1)}_{a_1b_1}
N^{(0)}_{a_2b_2} \textstyle \Bigl\{ n_1n_2 (n+3n_1) \binom{a-2}{a_1} +
n_2^2 (3n_1-n) \binom{a-2}{a_1-1} - 6 n_2^3 \binom{a-2}{a_1-2} \Bigr\}
\binom{b}{b_1} } \\
\shoveleft{{} + \tfrac{1}{12} {\displaystyle
\sum_{\substack{a=a_1+a_2\\b=b_1+b_2}}} N^{(0)}_{a_1b_1} N^{(0)}_{a_2b_2}
n_1 n_2^2 } \\
\textstyle \Bigl\{ n_1^2 (3-n_1) \binom{a-2}{a_1} + n_1n_2(n-3n_1-3)
\binom{a-2}{a_1-1} + n_2^2 (-n_1+n_2-6) \binom{a-2}{a_1-2} \Bigr\}
\binom{b}{b_1} \\
\shoveleft{{} + \tfrac{1}{2} \sum_{\substack{a=a_1+a_2+a_3\\b=b_1+b_2+b_3}}
\textstyle N^{(1)}_{a_1b_1} N^{(0)}_{a_2b_2} N^{(0)}_{a_3b_3} \Bigl\{
2n_1n_2^3n_3(n+3n_1-3n_2) \binom{a-2}{a_2,a_3-2} - 6 n_2^3n_3^3
\binom{a-2}{a_2,a_3} } \\
\textstyle {} + n_2^2n_3^2 (3n_1-n) \Bigl( n_1
\binom{a-2}{a_2-1,a_3-1} + n_2 \binom{a-2}{a_2,a_3-1} + n_3
\binom{a-2}{a_2-1,a_3} \Bigr) \Bigr\} \binom{b}{b_2,b_3} .
\end{multline*}
This relation determines the elliptic coefficient $N^{(1)}_{ab}$ for $a>0$
in terms of $N^{(1)}_{0,\frac{1}{2}a+b}$, the elliptic coefficients of
lower degree, and the rational coefficients.

To determine the coefficients $N^{(1)}_{0,b}$, we need the relation
obtained by evaluating Proposition \ref{Relation} on
$\om^2\boxtimes\om^2\boxtimes\om\boxtimes\om$: if $b\ge2$, then
\begin{multline*}
0 = N^{(1)}_{ab}  + \tfrac{1}{24} n(2n-1) N^{(0)}_{a+2,b-1}
+ \tfrac{1}{48} N^{(0)}_{a+4,b-2} \\
\shoveleft{ {} + \sum_{\substack{a+2=a_1+a_2\\b-1=b_1+b_2}} \textstyle
N^{(1)}_{a_1b_1} N^{(0)}_{a_2b_2} \textstyle \Bigl\{ n_2 \Bigl( n
\binom{a}{a_1} + n_2 \binom{a}{a_1-1} \Bigr) \binom{b-2}{b_1-1} } \\[-10pt]
\shoveright{ \textstyle {} - \frac{1}{6} \Bigl( n_1(6n_1-n_2)
\binom{a}{a_1} + n_2 (16n_1-n_2) \binom{a}{a_1-1} + 6n_2^2 \binom{a}{a_1-2}
\Bigr) \binom{b-2}{b_1} \Bigr\} } \\
\shoveleft{{} - \tfrac{1}{12} \sum_{\substack{a+4=a_1+a_2\\b-2=b_1+b_2}}
\textstyle
N^{(1)}_{a_1b_1} N^{(0)}_{a_2b_2} \Bigl( n_1 \binom{a}{a_1} + (2n_1-5n_2)
\binom{a}{a_1-1} + 6n_2 \binom{a}{a_1-2} \Bigr) \binom{b-2}{b_1} } \\
\shoveleft{ {} - \tfrac{1}{48} \sum_{\substack{a+4=a_1+a_2\\b-2=b_1+b_2}}
N^{(0)}_{a_1b_1} N^{(0)}_{a_2b_2} \textstyle \Bigl( n_1^3(n_1-1) \binom{a}{a_1}
+ n_1^2n_2(2n_1-2n_2+1) \binom{a}{a_1-1} } \\
\textstyle {} + n_1n_2^2(2n_1-2n_2+7) \binom{a}{a_1-2}
+ n_2^3(2n_1+5) \binom{a}{a_1-3} + n_2^4 \binom{a}{a_1-4} \Bigr)
\binom{b-2}{b_1} \\
\shoveleft{{} - \tfrac{1}{12}
\sum_{\substack{a+4=a_1+a_2+a_3\\b-2=b_1+b_2+b_3}}
\textstyle N^{(1)}_{a_1b_1} N^{(0)}_{a_2b_2} N^{(0)}_{a_3b_3}
\textstyle \Bigl\{ 3n_2n_3 \Bigl( n_2^2 \binom{a}{a_2,a_3-2} + n_3^2
\binom{a}{a_2-2,a_3} \Bigr) } \\
\shoveleft{ {} + \textstyle n_1 \Bigl( n_2^3 \binom{a}{a_2,a_3-4}
+ n_2^2(6n_1-n_3) \binom{a}{a_2-1,a_3-3} - 7n_2n_3^2 \binom{a}{a_2-2,a_3-2}
- 5n_3^3 \binom{a}{a_2-3,a_3-1} \Bigr) } \\
\shoveleft{\textstyle {} + \Bigl( n_2^3(n_1-5n_3) \binom{a}{a_2,a_3-3}
+ n_2^2n_3(5n_1-7n_3) \binom{a}{a_2-1,a_3-2} } \\
\textstyle {}+ n_2n_3^2 (5n_1-n_3) \binom{a}{a_2-2,a_3-1}
+ n_3^3(n_1+n_3) \binom{a}{a_2-3,a_3}\Bigr) \Bigr\} \binom{b-2}{b_2,b_3} .
\end{multline*}
This relation determine the coefficient $N^{(1)}_{0b}$ in terms of elliptic
coefficients of lower order and the rational coefficients, and thus
ultimately in terms of $N^{(0)}_{02}=1$, the number of lines between two
points.

Using these relation, we obtain the results of Table 2. Up to degree $3$,
Theorem A is easily seen to hold, since there are no elliptic space curves
of degrees $1$ and $2$, while all elliptic space curves of degree $3$ lie
in a plane.

It is well-known that there is one quartic elliptic space curve through $8$
general points, while the number of elliptic quartic space curves through
$16$ general lines was calculated by Vainsencher and Avritzer
(\cite{Vainsencher}; see also \cite{Avritzer}, which contains a correction
to \cite{Vainsencher}, bringing it into agreement with our calculation!).

\begin{table} \label{CP3}
\caption{Rational and elliptic Gromov-Witten invariants of $\CP^3$}
$$\begin{tabular}{|R|C|R|D{.}{}{2}|R|} \hline
n & (a,b) & N^{(0)}_{ab} & N^{(1)}_{ab} &
{\scriptstyle N^{(1)}_{ab}+(2n-1)N^{(0)}_{ab}/12} \\ \hline
1 & (0,2) & 1 & -.\frac{1}{12} & 0 \\
& (2,1) & 1 & -.\frac{1}{12} & 0 \\
& (4,0) & 2 & -.\frac{1}{6} & 0 \\[5pt]
2 & (0,4) & 0 & .0 & 0 \\
& (2,3) & 1 & -.\frac{1}{4} & 0 \\
& (4,2) & 4 & -1. & 0 \\
& (6,1) & 18 & -4.\frac{1}{2} & 0 \\
& (8,0) & 92 & -23. & 0 \\[5pt]
3 & (0,6) & 1 & -.\frac{5}{12} & 0 \\
& (2,5) & 5 & -2.\frac{1}{12} & 0 \\
& (4,4) & 30 & -12.\frac{1}{2} & 0 \\
& (6,3) & 190 & -78.\frac{1}{6} & 1 \\
& (8,2) & 1\,312 & -532.\frac{2}{3} & 14 \\
& (10,1) & 9\,864 & -3\,960. & 150 \\
& (12,0) & 80\,160 & -31\,900. & 1\,500 \\[5pt]
4 & (0,8) & 4 & -1.\frac{1}{3} & 1 \\
& (2,7) & 58 & -29.\frac{5}{6} & 4 \\
& (4,6) & 480 & -248. & 32 \\
& (6,5) & 4\,000 & -2\,023.\frac{1}{3} & 310 \\
& (8,4) & 35\,104 & -17\,257.\frac{1}{3} & 3\,220 \\
& (10,3) & 327\,888 & -156\,594. & 34\,674 \\
& (12,2) & 3259\,680 & -1\,515\,824. & 385\,656 \\
& (14,1) & 34\,382\,544 & -15\,620\,216. & 4\,436\,268 \\
& (16,0) & 383\,306\,880 & -170\,763\,640. & 52\,832\,040 \\[5pt]
5 & (0,10) & 105 & -36.\frac{3}{4} & 42 \\
& (2,9) & 1\,265 & -594.\frac{3}{4} & 354 \\
& (4,8) & 13\,354 & -6\,523.\frac{1}{2} & 3\,492 \\
& (6,7) & 139\,098 & -66\,274.\frac{1}{2} & 38\,049 \\
& (8,6) & 1\,492\,616 & -677\,808. & 441\,654 \\
& (10,5) & 16\,744\,080 & -7\,179\,606. & 5\,378\,454 \\
& (12,4) & 197\,240\,400 & -79\,637\,976. & 68\,292\,324 \\
& (14,3) & 2\,440\,235\,712 & -928\,521\,900. & 901\,654\,884 \\
& (16,2) & 31\,658\,432\,256 & -11\,385\,660\,384. & 12\,358\,163\,808 \\
& (18,1) & 429\,750\,191\,232 & -146\,713\,008\,096. & 175\,599\,635\,328 \\
& (20,0) & 6\,089\,786\,376\,960 & -1\,984\,020\,394\,752. &
2\,583\,319\,387\,968 \\ \hline
\end{tabular}$$
\end{table}

\end{document}